\documentclass[11pt]{article}
\usepackage{amssymb}
\begin{document}
\setlength{\oddsidemargin}{0cm}
\setlength{\evensidemargin}{0cm}
\baselineskip=20pt

\begin{center} {\Large\bf
Bijective 1-cocycles and Classification of 3-dimensional
Left-symmetric Algebras}
\end{center}

\bigskip

\begin{center}  { \large Chengming ${\rm Bai}^{1,2}$} \end{center}

\begin{center}{\it 1. Chern Institute of
Mathematics \& LMPC, Nankai University,  Tianjin 300071, P.R. China
}
\end{center}

\begin{center}{\it 2. Dept. of Mathematics, Rutgers, The State University of
New Jersey, Piscataway, NJ 08854, U.S.A.}\end{center}

\vspace{0.3cm}

\begin{center} {\large\bf   Abstract } \end{center}

Left-symmetric algebras have close relations with many important
fields in mathematics and mathematical physics. Their classification
is very complicated due to the nonassociativity. In this paper, we
re-study the correspondence between left-symmetric algebras and the
bijective 1-cocycles. Then a procedure is provided to classify
left-symmetric algebras in terms of classification of equivalent
classes of bijective 1-cocycles. As an example, the 3-dimensional
complex left-symmetric algebras are classified.

\vspace{0.2cm}



\section{Introduction}

Left-symmetric algebras (or under other names like Koszul-Vinberg
algebras, quasi-associative algebras, pre-Lie algebras and so on)
are a class of nonassociative algebras coming from the study of
several topics in geometry and algebra, such as rooted tree algebras
([C]), convex homogenous cones ([V]), affine manifolds and affine
structures on Lie groups ([Ko],[Ma]), deformation of associative
algebras ([G]) and so on. They are Lie-admissible algebras (in the
sense that the commutators define Lie algebra structures) whose left
multiplication operators form a Lie algebra.

Furthermore, left-symmetric algebras are a kind of natural algebraic
systems appearing in many fields in mathematics and mathematical
physics. Perhaps this is one of the most attractive and interesting
places. As it was pointed out in [CL], the left-symmetric algebra
``deserves more attention than it has been given". For example,
left-symmetric algebras appear as an underlying structure of those
Lie algebras that possess a phase space, thus ``they form a natural
category from the point of view of classical and quantum mechanics"
([Ku1-2]); they are the underlying algebraic structures of vertex
algebras ([BK]); there is a correspondence between left-symmetric
algebras and complex product structures on Lie algebras ([AS]),
which plays an important role in the theory of hypercomplex and
hypersymplectic manifolds ([Bar]); left-symmetric algebras have
close relations with certain integrable systems ([Bo],[LM]),
classical and quantum Yang-Baxter equation ([DM],[ESS],[GS],[Ku3]),
Poisson brackets and infinite-dimensional Lie algebras
([BN],[DN],[GD]), operads ([CL]), quantum field theory ([CK]) and so
on (see [Bu3] and the references therein).

On the other hand, it is hard to study left-symmetric algebras. Due
to the nonassociativity, there is neither a suitable representation
theory nor a complete structure theory like other classical algebras
such as associative algebras and Lie algebras. Even there exist
simple transitive left-symmetric algebras which combine the
simplicity and certain nilpotence ([H],[Bu1-2], or see the type
(${\rm D}_{-1}-10$) in section 3). In fact, many fundamental
problems have not been solved. Even the classification of complex
left-symmetric algebras is only up to dimension 2 ([BM1], [Bu2]).

Therefore we hope to get more examples which can be regarded as a
guide for further study. One of the ideas to get examples is to use
some well-known structures to obtain some left-symmetric algebras
(the so-called ``realization" theory). We have already obtained some
experiences. For example, a commutative associative algebra
$(A,\cdot)$ and its derivation $D$  can define a Novikov algebra
$(A,*)$ (which is a left-symmetric algebra with commutative right
multiplication operators) as follows ([GD],[BM4-5]):
$$x*y=x\cdot Dy,\;\;\forall x,y\in A.\eqno (1.1)$$

An analogue of the above construction in the version of Lie algebras
is related to the classical Yang-Baxter equation. In fact, a Lie
algebra $({\cal G},[\quad ])$ and a linear map $R:{\cal
G}\rightarrow {\cal G}$ satisfying the (operator form) of  classical
Yang-Baxter equation ([S])
$$[R(x),R(y)]=R([R(x),y]+[x,R(y)]), \;\;\forall x,\;y\in {\cal
G}\eqno (1.2)$$ can define a left-symmetric algebra $({\cal G},
*)$ as follows ([BM6],[GS],[Me]):
$$x*y=[R(x),y],\;\;\forall x,\;y\in {\cal G}.\eqno (1.3)$$
Moreover, equation (1.3) also gives an algebraic interpretation of
the so-called ``left-symmetry'': in some sense, the ``left-symmetry"
can be interpreted as a Lie bracket ``left-twisted'' by a classical
$r$-matrix. Furthermore, the above construction can be generalized
to any representation $(\rho, V)$ of a Lie algebra ${\cal G}$ as
follows. Let $T:V\rightarrow {\cal G}$ be a linear map satisfying
$$[T(u),T(v)]=T(\rho(T(u))v-\rho(T(v))u),\;\;\forall u,v \in {\cal
G}.\eqno (1.4)$$ Such a map is called an ${\cal O}$-operators in
[Ku3] which satisfies the (generalized) classical Yang-Baxter
equation (in fact, it is a solution of classical Yang-Baxter
equation on a larger Lie algebra ([Bai])). Then there exist
left-symmetric algebra structures in both $V$ and $T(V)\subset {\cal
G}$ given by
$$u*v=\rho(T(u))v;\;\;T(u)*T(v)=T(\rho(T(u))v),\;\;\forall u,v\in
V,\eqno (1.5)$$ respectively. These relations are not only useful
for the study of the left-symmetric algebras themselves such as
giving more examples as above and illuminating some interesting
properties, but also can provide those related topics with certain
algebraic and geometric interpretations.

Note that if the ${\cal O}$-operator $T$ appearing in equation (1.4)
is invertible, then the operator $T^{-1}$ is just a 1-cocycle
associated to the representation $(\rho, V)$ of ${\cal G}$. In fact,
there is a closer relation between left-symmetric algebras and
bijective 1-cocycles: there exists a compatible left-symmetric
algebra structure on a Lie algebra ${\cal G}$ if and only if ${\cal
G}$ has a bijective 1-cocycle. In this paper, we re-study the
correspondence between left-symmetric algebras and  bijective
1-cocycles.  Although most of the results have been already known
([Ki1], [Me]), our discussion can provide a procedure to classify
left-symmetric algebras using the representation theory of Lie
algebras. It is a ``linearization'' method which avoids classifying
the (non-linear) quadratic forms of structural constants.  In
particular, it is quite effective for the classification of complex
left-symmetric algebras in low dimensions, such as in dimension 3.

The paper is organized as follows. In Section 2, we re-study the
correspondence between the left-symmetric algebras and bijective
1-cocycles. In Section 3, we give the classification of
3-dimensional complex left-symmetric algebras.

Throughout this paper, without special saying, all algebras are of finite
dimension and over an algebraically closed field of characteristic 0.

\section{Left-symmetric algebras and bijective 1-cocycles}

\subsection{Preliminaries on left-symmetric algebras}

{\bf Definition 2.1}\quad Let $A$ be a vector space over a field ${\bf F}$
equipped with a bilinear product $(x,y)\rightarrow xy$. $A$ is called a
left-symmetric algebra if for any $x,y,z\in A$, the associator
$$(x,y,z)=(xy)z-x(yz)\eqno (2.1)$$
is symmetric in $x,y$, that is,
$$(x,y,z)=(y,x,z),\;\;{\rm or}\;\;{\rm
equivalently}\;\;(xy)z-x(yz)=(yx)z-y(xz).\eqno (2.2)$$

For a left-symmetric algebra $A$, the commutator
$$[x,y]=xy-yx,\eqno (2.3)$$
defines a Lie algebra ${\cal G}={\cal G}(A)$, which is called the
sub-adjacent Lie algebra of $A$. For any $x,y\in A$, let $L_x$ and
$R_x$ denote the left and right multiplication operator
respectively, that is, $L_x(y)=xy,\;R_x(y)=yx$. Then the
left-symmetry (2.2) is just
$$[L_x,L_y]=L_{[x,y]},\;\;\forall x,y\in A, \eqno (2.4)$$
which means that $L:{\cal G}(A)\rightarrow
gl({\cal G}(A))$ with $x\rightarrow L_x$
gives a (regular) representation of the Lie algebra ${\cal G}(A)$.

Some subclasses of left-symmetric algebras are very important.

{\bf Definition 2.2} \quad Let $A$ be a left-symmetric algebra.

(1) If $A$ has no ideals except itself and zero, then $A$ is called simple.
$A$ is called semisimple if $A$ is the direct sum of
simple left-symmetric algebras.

(2) If for every $x\in A$, $R_x$ is nilpotent, then $A$ is said to
be transitive or complete. The transitivity corresponds to the
completeness of an affine manifold ([Ki1],[Me]). Moreover, the
sub-adjacent Lie algebra of a transitive left-symmetric algebra is
solvable.

(3) If for every $x,y\in A$, $R_xR_y=R_yR_x$, then $A$ is called a
Novikov algebra. It was introduced in connection with the Poisson
brackets of hydrodynamic type and Hamiltonian operators in the
formal variational calculus ([BN],[BM3-6],[GD],[O],[X],[Z]).

(4) If for every $x,y,z\in A$, the associator $(x,y,z)$ is
right-symmetric, that is, $(x,y,z)=(x,z,y)$, then $A$ is said to be
bi-symmetric. It is just the assosymmetric ring in the study of near
associative algebras ([Kl],[BM2]).

\subsection{Bijective 1-cocycles}

{\bf Definition 2.3}\quad Let ${\cal G}$ be a Lie algebra and
$f:{\cal G}\rightarrow gl(V)$ be a representation of ${\cal G}$. A
1-cocycle $q$ associated to $f$ is defined as a linear map from
${\cal G}$ to $V$ satisfying
$$q[x,y]=f(x)q(y)-f(y)q(x),\forall x,y\in {\cal G}.\eqno (2.5)$$
We denote it by $(f,q)$. In addition, if $q$ is a linear isomorphism
(thus $\dim  V=\dim{\cal G}$), $(f,q)$ is said to be bijective.

Let $(f,q)$ be a bijective 1-cocycle, then it is easy to see that
$$x*y=q^{-1}(f(x)q(y)),\;\;\forall x, y\in {\cal G}.\eqno (2.6)$$
defines a left-symmetric algebra on ${\cal G}$ ([Me]). Conversely,
for a left-symmetric algebra ${\cal G}$, $(L, id)$ is a bijective
1-cocycle of ${\cal G}$. Hence we have the following maps:

$\Phi: {\bf A}=\{{\rm bijective}\;\;{1-}{\rm
cocycles}\}\longrightarrow {\bf B} =\{\;{\rm left}-{\rm symmetric
}\;\;{\rm algebras}\}$

$\Psi:{\bf B}=\{\;{\rm left}-{\rm symmetric }\;\;{\rm
algebras}\}\longrightarrow {\bf A}=\{{\rm bijective}\;\;{1-}{\rm
cocycles}\}$

{\bf Definition 2.4}\quad Let ${\cal G}$ be a Lie algebra. Let
$(f_1, V_1)$ and $(f_2, V_2)$ be two linear representations and
$q_1, q_2$ be bijective 1-cocycles associated to $f_1, f_2$
respectively. $(f_1, V_1)$ is isomorphic ($\cong$) to $(f_2,V_2)$ if
there exists a linear isomorphism $g:V_1\rightarrow V_2$ such that
$f_2=g f_1 g^{-1}$. We call them equivalent ($\sim$) if there exists
an automorphism $T$ of ${\cal G}$ such that $(f_1T, V_1)\cong (f_2,
V_2)$. $(f_1,q_1)$  is isomorphic ($\cong$) to $(f_2,q_2)$ if there
exists a linear isomorphism $g:V_1\rightarrow V_2$ such that $f_2=g
f_1 g^{-1}$ and $q_2=gq_1$. We call them equivalent ($\sim$) if
there exists an automorphism $T$ of ${\cal G}$ such that $(f_1T,
q_1T)\cong (\rho_2, q_2)$, that is, there exist a linear isomorphism
$g:V_1\rightarrow V_2$ and an automorphism $T:{\cal G}\rightarrow
{\cal G}$ such that $f_2=gf_1Tg^{-1}$ and $q_2=gq_1T$.

On the other hand, recall that two left-symmetric algebras $({\cal
G}_1,*)$ and $({\cal G}_2,\cdot)$ are isomorphic (denoted by
$\cong$) if there exists a linear isomorphism $F:{\cal
G}_1\rightarrow {\cal G}_2$ such that $F(x*y)=F(x)\cdot F(y)$ for
any $x,y\in {\cal G}_1$.

{\bf Theorem 2.1}\quad The maps $\Phi$ and $\Psi$ induce a bijection
between the set of the isomorphism classes of bijective 1-cocycles
of ${\cal G}$ and the set of left-symmetric algebras on ${\cal G}$.
Under this correspondence equivalent bijective 1-cocycles are mapped
to isomorphic left-symmetric algebras and vice versa. That is,
$${\bf A}/{\cong}\longleftrightarrow {\bf B};\;\;{\bf A}/{\sim}\longleftrightarrow {\bf
B}/{\cong}.\eqno (2.7)$$

{\bf Proof}\quad Let $(f_1,q_1)$ and $(f_2, q_2)$ be two isomorphic
bijective 1-cocycles. Then there exists a linear isomorphism $g$
such that $gf_1=f_2g$ and $q_2=gq_1$. We can know their
corresponding left-symmetric algebras coincide since
$$x*y=q_1^{-1}(f_1(x)q_1(y))=(q_2^{-1}g)[(g^{-1}f_2(x)g)(g^{-1}q_2)(y)]=q_2^{-1}(f_2(x)q_2(y)),\;\;\forall x,y\in {\cal G}.$$
Therefore the map $\Phi$ is defined on the set of isomorphism
classes of bijective 1-cocycles.

Let $(f,q)$ be a bijective 1-cocycle, then $\Psi\Phi(f,q)\cong
(f,q)$ by $g=q$. Conversely, from the definitions, we know that
$\Phi\Psi (*)=*$, that is, $\Phi\Psi$ maps any left-symmetric
algebra to itself. Hence the correspondence is proved.

Now we prove that equivalent bijective 1-cocycles correspond to
isomorphic left-symmetric algebras. Let $(f_1,q_1)$ and $(f_2, q_2)$
be two equivalent bijective 1-cocycles. Then there exists a linear
isomorphism there exist a linear isomorphism $g:V_1\rightarrow V_2$
and an automorphism $T:{\cal G}\rightarrow {\cal G}$ such that
$f_2=gf_1Tg^{-1}$ and $q_2=gq_1T$. Their corresponding
left-symmetric algebra are given by
$$x*_1y=q_1^{-1}(f_1(x)q_1(y)),\;\;\forall x,y\in {\cal G}$$
and
$$x*_2y=q_2^{-1}(f_2(x)q_2(y))=(T^{-1}q_1^{-1}g^{-1})[(gf_1T(x)g^{-1})(gq_1T(y))],\;\;\forall
x,y\in {\cal G},$$ respectively. So we have
$$Tx*_1Ty=T(x*_2y).\;\;\forall x,y\in {\cal G}.$$
Hence these two left-symmetric algebras are isomorphic by $T$.
Conversely, let $F$ be an (left-symmetric algebra) automorphism of
${\cal G}$. Obviously, $F$ is also a Lie algebra automorphism of
${\cal G}$ and $(LF, F)$ is a bijective 1-cocycle of ${\cal G}$
corresponding to the image of $F$. Then the bijective 1-cocycle $(L,
id)$ is equivalent to the bijective 1-cocycle $(LF,F)$ by $g=id$.
\hfill $\Box$

Hence the classification of left-symmetric algebras in the sense of
isomorphism is as the same as the classification of bijective
1-cocycles in the sense of equivalence.

{\bf Remark 1}\quad The above correspondence is similar to the
correspondence between left-symmetric algebras and \'etale affine
representations given in [Bau].

{\bf Remark 2} As in the introduction, when $(f, q)$ is a bijective
1-cocycle of a Lie algebra ${\cal G}$, then $q^{-1}$ is an ${\cal
O}$-operator associated to $f$, that is, $q^{-1}$ satisfies the
(generalized) classical Yang-Baxter equation. Moreover, equation
(2.6) coincides with the latter part of equation (1.5) since
$$x*y=q^{-1}(u)*q^{-1}(v)=q^{-1}(f(q^{-1}(u))v)=q^{-1}(f(x)q(y)),\eqno (2.8)$$
where $x=q^{-1}(u), y=q^{-1}(v)$.

\subsection{Classification Problems}

Due to the nonassociativity, it is very difficult to classify
left-symmetric algebras. A natural way is to classify the structural
constants, which has been used in dimension 2 ([BM1] or [Bu2]).
However, it can not be extended to the higher dimensions since it
involves the classification of quadratic forms of structural
constants, which is very complicated due to nonlinearity, even in
dimension 3. Moreover, unlike associative algebras or Lie algebras,
there is not a complete structure theory. For example, although
there are several definitions of radicals ([Bu1-2],[H],[Me]), none
is good enough. In fact, up to now, there are only 2-dimensional
complex left-symmetric algebras and some special cases in higher
dimensions (for example, transitive cases on nilpotent Lie algebras
up to dimension 4 ([Ki1-2]) and bi-symmetric algebras ([BM2]) and
Novikov algebras ([BM3]) up to dimension 3) have been classified.

From the relation between left-symmetric algebras and bijective
1-cocycles, we can solve this problem by classifying the
equivalent bijective 1-cocycles. In fact, we can divide the classification into
several steps:

Step 1: Classify Lie algebras. This has been done in certain
low dimensions and some special cases ([J],[SW]).

Step 2: Let ${\cal G}$ be a given Lie algebra with a basis $\{
e_1,\cdots,
e_n\}$. Compute its automorphism group ${\rm Aut}({\cal G})$.
For a representation $f:{\cal
G}\rightarrow gl(V)$ ($\dim V=\dim
{\cal G}$) with a basis $\{ v_1,\cdots, v_n\}$ of $V$,  we can let
$f(x)=(f_{ij}(x))$ for any $x\in {\cal G}$, where $f_{ij}:{\cal
G}\rightarrow {\bf F}$ be
linear functions. On the other hand, let $q:{\cal G}\rightarrow V$ be
a 1-cocycle, then we
can  let $q(x)=\sum\limits_{k=1}^n
A_k(x)v_k$, where $A_k:{\cal G}\rightarrow {\bf F}$ are linear functions.
The conditions of the representation $f$ and the 1-cocycle $q$ can give
a series of equations for linear functions $f_{ij}$ and $A_k$.

Step 3: Classify the linear functions $f_{ij}$ under the sense of
equivalence through the basis transformations of
$V$ and the basis transformation of ${\cal G}$ which is compatible with the
automorphism group of ${\cal G}$.

Step 4: For a given representation obtained in step 3, find all
the corresponding bijective 1-cocycles (that is, the determinant of $(A_j(e_i))$
is non-zero).

Step 5: Classify those bijective 1-cocycles in step 4 and their
corresponding
left-symmetric algebras.

Although it seems that it is more complicated to classify
bijective 1-cocycles than left-symmetric algebras themselves, in
fact, there are certain advantages: every above step only involves
linear equations, thus avoiding the classification of the
nonlinearity of structural constants; the whole classification is
like a kind of ``variable separated'' (in particular in step 3 and
step 5). The whole process is like a kind of ``linearization'' of
classifying structural constants of left-symmetric algebras.
Moreover, this method  can be extended to use the extensions of
left-symmetric algebras ([Ki2]).

The above procedure will be quite effective to classify some
left-symmetric algebras over the complex number ${\bf C}$. As an
example, we give the classification of 3-dimensional complex
left-symmetric algebras in Section 3.

\section{The classification of 3-dimensional complex left-symmetric
algebras}

It is well-known that there does not exist any left-symmetric
algebra structure on a complex semisimple Lie algebra (cf. [Me]).
Hence, over the complex field ${\bf C}$, besides 3-dimensional
simple Lie algebra $sl(2,{\bf C})$, up to isomorphisms, there are
the following (non-isomorphic) Lie algebras ([J]): (we only give the
non-zero products)

(a) Abelian Lie algebra;

(b) Heisenberg Lie algebra ${\cal H}=< e_1,e_2,e_3|[e_1,e_2]=e_3>$;

(c) ${\cal N}=< e_1, e_2, e_3|[e_3,e_2]=e_2>$, which is a
direct sum of 2-dimensional non-abelian Lie algebra and 1-dimensional
center;

(d) ${\cal D}_l=< e_1, e_2,e_3|[e_3,e_1]=e_1,[e_3,e_2]=le_2>,\;
0<|l|<1$ or $l=e^{i\theta}, 0\leq \theta\leq \pi$;

(e) ${\cal E}=<e_1, e_2, e_3|[e_3,e_1]=e_1, [e_3, e_2]=e_1+e_2>$.

All of the above Lie algebras are solvable. Let ${\cal G}$ be one of
these algebras and $f:{\cal G}\rightarrow gl(V)$ be a
representation. Thus, according to Lie Theorem ([J]), there exists a
basis $\{v_1,v_2,v_3\}$ of $V$ such that
$$f(x)=\left(\matrix{ f_{11}(x) & 0 &0\cr
f_{21}(x) & f_{22} (x) &0\cr f_{31}(x) & f_{32}(x)
&f_{33}(x)\cr}\right),\forall x\in {\cal G}\eqno (3.1)$$ where
$f_{11},f_{22},f_{33},f_{21},f_{32},f_{31}$ are linear functions of
${\cal G}$. Let $q:{\cal G}\rightarrow V$ be 1-cocycle:
$$q(x)=A_1(x)v_1+A_2(x)v_2+A_3(x)v_3,\eqno (3.2)$$
where $A_1, A_2,A_3$ are linear functions of ${\cal G}$. The matrix
associated to $q$ is
defined as $C=(A_j(e_i))=\left(\matrix{ A_1(e_1) & A_2(e_1) & A_3(e_1)\cr
A_2(e_1) & A_2(e_2) & A_3(e_3)\cr A_3(e_1) & A_3(e_2) &
A_3(e_3)\cr}\right)$. $q$ is bijective if and only if $\det C\ne 0$.

For a left-symmetric algebra ${\cal G}$, the (form)  characteristic
matrix is defined as
$$M=\left(\matrix{ \sum_{k=1}^n a_{11}^ke_k &\cdots & \sum_{k=1}^n
a_{1n}^ke_k\cr \cdots &\cdots &\cdots\cr \sum_{k=1}^n a_{n1}^ke_k &
\cdots &\sum_{k=1}^n a_{nn}^ke_k\cr}\right), \eqno (3.3)$$ where
$\{e_i\}$ is a basis of ${\cal G}$ and $e_ie_j=\sum_{k=1}^n
a_{ij}^ke_k$.

The left-symmetric algebras on abelian Lie algebras are commutative
associative algebras which are classified in dimension 3 in [BM3].
In the next subsections, we give the classification of 3-dimensional
complex left-symmetric algebras on the Lie algebras (b)-(e)
according to the procedure given in last section. As an explanation,
we give a detailed and explicit demonstration for the left-symmetric
algebras on ${\cal H}$, whereas we omit the length proof for other
cases since the proof is similar.

\subsection{The left-symmetric algebras on ${\cal H}$}

The automorphism group of ${\cal H}$ is
$${\rm Aut}({\cal H})=\{
\left( \matrix{ a_{11} & a_{12} & a_{13}\cr  a_{21} & a_{22} &
a_{23}\cr 0 & 0&
a_{11}a_{22}-a_{12}a_{21}\cr}\right)|a_{11}a_{22}-a_{12}a_{21}\ne
0\}.\eqno (3.4)$$ It is easy to show that $f: {\cal H}\rightarrow
gl(V)$ defined by equation (3.1) is a representation if and only if
it satisfies the following conditions:
\begin{eqnarray*}
&&f_{11}(e_3)=f_{22}
(e_3)=f_{33}(e_3)=f_{21}(e_3)=f_{32}(e_3)=0;\\
&&f_{21}(e_1)(f_{22}-f_{11})(e_2)+f_{21}(e_2)(f_{11}-f_{22})(e_1)=0;\\
&&f_{32}(e_1)(f_{33}-f_{22})(e_2)+f_{32}(e_2)(f_{22}-f_{33})(e_1)=0;\\
&&f_{31} (e_3)(f_{11}-f_{33})(e_1)=f_{31} (e_3)
(f_{11}-f_{33})(e_2)=0;\\
&&f_{31}(e_3)=f_{31}(e_1)
(f_{33}-f_{11})(e_2)+f_{31}(e_2)(f_{11}-f_{33})(e_1)+f_{21}(e_1)f_{32}(e_2)-f_{21}(e_2)f_{32}(e_1).
\end{eqnarray*}
The 1-cocycle $C=(A_j(e_i))$ satisfies the following conditions:
\begin{eqnarray*}
&&A_3(e_3)=0,\;-A_3(e_1)f_{33}(e_2)+A_3(e_2)f_{33}(e_1)=0;\\
&&A_2(e_3)f_{22}(e_1)=A_2(e_3)f_{22}(e_2)=0;\\
&&A_2(e_1)f_{22}(e_2)-A_2(e_2)
f_{22}(e_1)+A_3(e_1)f_{32}(e_2)-A_3(e_2)f_{32}(e_1)=-A_2(e_3);\\
&&-A_1(e_3)f_{11}(e_1)-A_2(e_3)f_{21}(e_1)+A_3(e_1)f_{31}(e_3)=0;\\
&&-A_1(e_3)f_{11}(e_2)-A_2(e_3)f_{21}(e_2)+A_3(e_2)f_{31}(e_3)=0;\\
&&A_1(e_1)f_{11}(e_2)-A_1(e_2)f_{11}(e_2)+A_2(e_1)f_{21}(e_2)-A_2(e_2)f_{21}
(e_1)+A_3(e_1)f_{31}(e_2)-A_3(e_2)f_{31}(e_1)\\
&&\mbox{}\hspace{0.5cm} =-A_1(e_3).
\end{eqnarray*}

{\bf Proposition 3.1}\quad If $f_{31} (e_3)\ne 0$, then the
representation must be equivalent to one of the following cases:
{\small $${\rm (AI)} \;f(e_1)=\left(\matrix{f_{11} (e_1) & 0 &0\cr 1
&f_{11}(e_1) &0\cr 0 &
0&f_{11}(e_1)\cr}\right),f(e_2)=\left(\matrix{f_{11} (e_2) & 0 &0\cr
0&f_{11}(e_2) &0\cr 0 &
1&f_{11}(e_2)\cr}\right),f(e_3)=\left(\matrix{ 0& 0 &0\cr 0 &0 &0\cr
1 & 0&0\cr}\right).$$}

{\bf Proof}\quad Since $f_{31}(e_3)\ne 0$, we can let
$f_{31}(e_3)=1$ through
$$v_1\rightarrow f_{31}(e_3)v_1,\;v_2\rightarrow v_2,\;v_3\rightarrow v_3.$$
Thus, by the equations of $f$, we know that $f_{11}=f_{33}$ and
$f_{32}(e_2)f_{21}(e_1) -f_{32}(e_1)f_{21}(e_2)=1$. We claim that
$f_{11}=f_{22}$. Otherwise, we can suppose $f_{11}(e_1)\ne
f_{22}(e_1)$. Then through
$$v_1\rightarrow v_1,\; v_2\rightarrow
v_2-\frac{f_{21}(e_1)}{(f_{11}-f_{22})(e_1)}v_1,\;\;v_3 \rightarrow
v_3-\frac{f_{32}(e_1)}{(f_{22}-f_{11})(e_1)}v_2$$ we can let
$f_{21}(e_1)=f_{32}(e_1)=0$, which is a contradiction. On the other
hand, we can let $f_{31}(e_1)=f_{31}(e_2)=0$ through the following
transformation which is in ${\rm Aut}({\cal H})$:
$$e_1\rightarrow e_1-f_{31}(e_1)e_3,\;e_2\rightarrow e_2-f_{31}(e_2)e_3,
e_3\rightarrow e_3.$$ Since $f_{21}$ and $f_{32}$ can not be zero,
without losing generality, we suppose $f_{21}(e_1)\ne 0$. Thus by
$$e_1\rightarrow e_1,\; e_2\rightarrow
e_2-\frac{f_{21}(e_2)}{f_{21}(e_1)}e_1,\;e_3\rightarrow e_3,$$ we
can let $f_{21}(e_2)=0$. Then $f_{21}(e_1)f_{32}(e_2)=1$. By
$$v_1\rightarrow v_1,\;v_2\rightarrow f_{32}(e_2)v_2,\;v_3\rightarrow v_3,$$
we can let $f_{32}(e_2)=f_{21}(e_1)=1$. Finally, we can get (AI) by
$$e_1\rightarrow e_1-f_{32}(e_1)e_2,\; e_2\rightarrow e_2,\;e_3\rightarrow
e_3.\eqno \Box$$

{\bf Proposition 3.2}\quad For a representation given in the case
(AI), there exist bijective 1-cocycles if and only if it is
equivalent to one of the following cases:

(AI-1) $f_{11}(e_1)=1,f_{11}(e_2)=0$. There is only one bijective
1-cocycle up to equivalence
$$C=\left(\matrix{0 &0 & 1\cr 0 & 1 &0\cr 1&0 &0\cr}\right)\Longrightarrow
{\rm (H-1)}\;\;\left(\matrix{ e_1 &e_2+e_3 & e_3\cr e_2 &0 &0\cr
e_3& 0 &0\cr}\right).$$

(AI-2) $f_{11}(e_1)=f_{11}(e_2)=1$. There is only one bijective
1-cocycle up to equivalence
$$C=\left(\matrix{0 &0 & 1\cr 0 & 1 &1\cr 1&0 &0\cr}\right)\Longrightarrow
{\rm (H-2)'}\;\left(\matrix{ e_1 &e_2+e_3 & e_3\cr e_2& 2e_2-e_1
&e_3 \cr e_3& e_3 &0\cr}\right)\cong {\rm (H-2)}\;\left(\matrix{ e_1
&e_2+e_3 & e_3\cr e_2 &e_3 &0\cr e_3& 0 &0\cr}\right).$$

{\bf Proof}\quad For a representation given in (AI), the equations
for $C=(A_j(e_i))$ reduce to the following equations:
\begin{eqnarray*}
&&A_3(e_3)=A_2(e_3)=0;f_{11}(e_1)A_3(e_2)=f_{11}(e_2)A_3(e_1); \\
&&A_3(e_1)=A_2(e_2)f_{11}(e_1)-A_2(e_1)f_{11}(e_2)=A_1(e_3)f_{11}(e_1);\\
&&A_3(e_2)=A_1(e_3)f_{11}(e_2);
A_1(e_3)=A_2(e_2)+A_1(e_2)f_{11}(e_1)-A_1(e_1)f_{11}(e_2).
\end{eqnarray*}
\indent If $f_{11}(e_1)=0$, then $A_3(e_1)=0$. If $f_{11}(e_2)\ne
0$, then $A_2(e_1)=0$, which leads to $\det C=0$. If
$f_{11}(e_2)=0$, then $A_3(e_2)=0$, which also leads to $\det C=0$.

If $f_{11}(e_1)\ne 0$, then we can let $f_{11}(e_1)=1$ by
$$e_1\rightarrow \frac{1}{f_{11}(e_1)}e_1, e_2\rightarrow e_2,
e_3\rightarrow \frac{1}{f_{11}(e_1)}e_3;\;v_1\rightarrow
\frac{1}{f_{11}(e_1)}v_1, v_2\rightarrow v_2,v_3\rightarrow v_3.$$
If $f_{11}(e_2)=0$, this is just the case (AI-1). At the time, the
corresponding bijective 1-cocycles are given by
$C=\left(\matrix{A_1(e_1) & A_2(e_1) & A_3(e_1)\cr 0 & A_3(e_1)
&0\cr A_3(e_1) & 0 &0\cr}\right)$ with $ A_3(e_1)\ne 0$. The
corresponding left-symmetric algebras are $\left(\matrix{
e_1+\frac{A_2(e_1)}{A_3(e_1)}e_3 &e_2+e_3 & e_3\cr e_2 &0 &0\cr e_3&
0 &0\cr}\right)$. However, they are isomorphic to (H-1) through
$$e_1\rightarrow e_1-\frac{A_2(e_1)}{A_3(e_1)}e_3, e_2\rightarrow e_2,
e_3\rightarrow e_3,$$ which can be given by $C=\left(\matrix{0 &0 &
1\cr 0 & 1 &0\cr 1&0 &0\cr}\right)$ directly.

\noindent If $f_{11}(e_2)\ne 0$, then we can let $f_{11}(e_2)=1$ by
$$ e_1\rightarrow e_1, e_2\rightarrow \frac{1}{f_{11}(e_2)}e_2,
e_3\rightarrow \frac{1}{f_{11}(e_2)}e_3;\; v_1\rightarrow v_1,
v_2\rightarrow v_2, v_3\rightarrow f_{11}(e_2)v_3,$$ which is the
case (AI-2). The corresponding bijective 1-cocycles are given by
$$C=\left(\matrix{A_1(e_1) & A_2(e_1) & A_3(e_1)\cr A_1(e_1)-A_2(e_1) &
A_2(e_1)+A_3(e_1) &A_3(e_1)\cr A_3(e_1) & 0 &0\cr}\right)$$ with
$A_3(e_1)\ne 0$. The corresponding left-symmetric algebras are
$$\left(\matrix{
e_1+\frac{A_2(e_1)}{A_3(e_1)}e_3
&e_2+\frac{A_2(e_1)+A_3(e_1)}{A_3(e_1)}e_3 & e_3\cr
e_2+\frac{A_2(e_1)}{A_3(e_1)}e_3
&2e_2-e_1+\frac{A_2(e_1)}{A_3(e_1)}e_3 &e_3\cr e_3& e_3
&0\cr}\right).$$ However, they are isomorphic to (H-2)' through
$$e_1\rightarrow e_1-\frac{A_2(e_1)}{A_3(e_1)}e_3, e_2\rightarrow e_2,
e_3\rightarrow e_3,$$ which can be given by $C=\left(\matrix{0 &0 &
1\cr 0 & 1 &1\cr 1&0 &0\cr}\right)$ directly. Notice that (H-2)'
$\cong$ (H-2) through
$$e_1\rightarrow e_1,\;e_2\rightarrow e_1-e_2, e_3\rightarrow
-e_3.\eqno\Box$$

{\bf Proposition 3.3}\quad If $f_{31}(e_3)=0$, then the equivalent
classes of the representations of ${\cal H}$ are divided into the
following cases (maybe there are equivalent classes in the same
case, but not equivalent in different cases):

${\rm (BI)}\; f(e_1)=\left(\matrix{f_{11} (e_1) & 0 &0\cr 0
&f_{22}(e_1) &0\cr 0 &
0&f_{33}(e_1)\cr}\right),f(e_2)=\left(\matrix{f_{11} (e_2) & 0 &0\cr
0 &f_{22}(e_2) &0\cr 0 & 0&f_{33}(e_2)\cr}\right),f(e_3)=0.$

${\rm (BII)}\; f(e_1)=\left(\matrix{f_{11} (e_1) & 0 &0\cr 0
&f_{11}(e_1) &0\cr 0 &
0&f_{33}(e_1)\cr}\right),f(e_2)=\left(\matrix{f_{11} (e_2) & 0 &0\cr
1&f_{11}(e_2) &0\cr 0 & 0&f_{33}(e_2)\cr}\right),f(e_3)=0.$

${\rm (BIII)}\;f(e_1)=\left(\matrix{f_{11} (e_1) & 0 &0\cr 0
&f_{11}(e_1) &0\cr 0 &
0&f_{11}(e_1)\cr}\right),f(e_2)=\left(\matrix{f_{11} (e_2) & 0 &0\cr
1 &f_{11}(e_2) &0\cr 0 & 1&f_{11}(e_2)\cr}\right),f(e_3)=0.$

${\rm (BIV)}\;f(e_1)=\left(\matrix{f_{11} (e_1) & 0 &0\cr 1
&f_{11}(e_1) &0\cr 0 &
0&f_{11}(e_1)\cr}\right),f(e_2)=\left(\matrix{f_{11} (e_2) & 0 &0\cr
0 &f_{11}(e_2) &0\cr 1 & 0&f_{11}(e_2)\cr}\right),f(e_3)=0.$

${\rm (BV)}\;f(e_1)=\left(\matrix{f_{11} (e_1) & 0 &0\cr 0
&f_{11}(e_1) &0\cr 0 &
1&f_{11}(e_1)\cr}\right),f(e_2)=\left(\matrix{f_{11} (e_2) & 0 &0\cr
0 &f_{11}(e_2) &0\cr 1 & 0&f_{11}(e_2)\cr}\right),f(e_3)=0.$

${\rm (BVI)}\;f(e_1)=\left(\matrix{f_{11} (e_1) & 0 &0\cr 0
&f_{11}(e_1) &0\cr 1 &
0&f_{11}(e_1)\cr}\right),f(e_2)=\left(\matrix{f_{11} (e_2) & 0 &0\cr
1 &f_{11}(e_2) &0\cr 0 & 1&f_{11}(e_2)\cr}\right),f(e_3)=0.$

{\bf Proof}\quad We give the sketch of proof here. The detailed
discussion is as the same as the discussion in the case (AI). Since
$f(e_3)=0$, we only need to consider $f(e_1)$ and $f(e_2)$.
Moreover, there is a kind of symmetry between $f(e_1)$ and $f(e_2)$
since $e_1\rightarrow e_2, e_2\rightarrow e_1, e_3\rightarrow -e_3$
is in ${\rm Aut}({\cal H})$. Hence we only need to consider the
Jordan canonical forms of $f(e_1)$.

$f(e_1)$ is diagonalized. Then the equations of $f$  reduce to
$$f_{21}(e_2)(f_{11}-f_{22})(e_1)=f_{32}(e_2)(f_{22}-f_{33})(e_1)=f_{31}(e_2)(f_{11}-f_{33})=0.$$
Thus we can consider the Jordan canonical form of $f(e_2)$. First of
all, $f(e_2)$ is also diagonalized. This is the case (BI). Secondly,
$f(e_2)$ is the type $\left(\matrix{f_{11} (e_2) & 0 &0\cr
1&f_{11}(e_2) &0\cr 0 & 0&f_{33}(e_2)\cr}\right)$, then
$f_{11}(e_1)=f_{22}(e_1)$, which is the case (BII). For the other
positions of Jordon blocks of $f(e_2)$, it is easy to show that they
are isomorphic. Finally, $f(e_2)$ is the type $\left(\matrix{f_{11}
(e_2) & 0 &0\cr 1&f_{11}(e_2) &0\cr 0 & 1&f_{33}(e_2)\cr}\right)$,
then $f_{11}(e_1)=f_{22}(e_1)=f_{33}(e_1)$. This is the case (BIII).

$f(e_1)$ is the type $\left(\matrix{f_{11} (e_1) & 0 &0\cr
1&f_{11}(e_1) &0\cr 0 & 0&f_{33}(e_1)\cr}\right)$. Then from the
equations of $f$, we can get
$$f_{32}(e_2)=0,\;\; f_{31}(e_2)(f_{11}-f_{33})(e_1)=0.$$
We can let $f_{21}(e_2)=0$ through
$$e_1\rightarrow e_1, e_2\rightarrow e_2-f_{21}(e_2) e_1, e_3\rightarrow
e_3.$$ If $f_{31}(e_2)=0$, by symmetry of $e_1, e_2$, it is
equivalent to the case (BII). If $f_{31}(e_2)\ne 0$, then
$f_{11}(e_1)=f_{33}(e_1)$ and moreover, we can let $f_{31}(e_2)=1$
by
$$e_1\rightarrow e_1, e_2\rightarrow \frac{1}{f_{31}(e_2)} e_2,
e_3\rightarrow \frac{1}{f_{31}(e_2)} e_3,$$ which is just the case
(BIV). Similarly, for other positions of the Jordan blocks of
$f(e_1)$ with the same type such as $f(e_1)=\left(\matrix{f_{11}
(e_1) & 0 &0\cr 0&f_{22}(e_1) &0\cr 0 & 1&f_{22}(e_1)\cr}\right)$
and $f(e_1)=\left(\matrix{f_{11} (e_1) & 0 &0\cr 0&f_{22}(e_1) &0\cr
1 & 0&f_{11}(e_1)\cr}\right)$, we can get the case (BV) and (BVI)
respectively.

$f(e_1)$ is the type $\left(\matrix{f_{11} (e_1) & 0 &0\cr
1&f_{11}(e_1) &0\cr 0 & 1&f_{11}(e_1)\cr}\right)$. Then from the
equations of $f$, we know
$$f_{11}(e_2)=f_{22}(e_2)=f_{33}(e_2), f_{32}(e_2)=f_{21}(e_2).$$
We can let $f_{32}(e_2)=f_{21}(e_2)=0$ by
$$e_1\rightarrow e_1,e_2\rightarrow e_2-f_{32}(e_2)e_1, e_3\rightarrow
e_3.$$ Then it is equivalent to the case (BIII) if $f_{31}(e_2)=0$
or  the case (BVI) if $f_{31}(e_2)\ne 0$.\hfill $\Box$

{\bf Proposition 3.4}\quad For a representation of ${\cal H}$ given
in the above cases respectively, there exist bijective 1-cocycles if
and only if it is equivalent to one of the following corresponding
cases: (we give the classification of 1-cocycles up to equivalence
and the corresponding left-symmetric algebras respectively)

\noindent Case (BI): there does not exist any bijective 1-cocycle;

\noindent Case (BII):\quad (BII-1) $f_{11}(e_1)=0,
f_{33}(e_1)=1,f_{11}(e_2)=f_{33}(e_2)=0$.
$$C=\left(\matrix{0 &-1 & 1\cr 0 & 1 &0\cr 1&0 &0\cr}\right)\Longrightarrow
{\rm (H-3)'}\;\left(\matrix{ e_1+e_2 &0 & 0\cr -e_3 &e_3 &0\cr 0& 0
&0\cr}\right)\cong {\rm (H-3)}\;\left(\matrix{ e_1 &e_3 & 0\cr 0
&e_3 &0\cr 0& 0 &0\cr}\right).$$ \indent (BII-2)
$f_{11}(e_1)=f_{33}(e_1)=0,f_{11}(e_2)=0, f_{33}(e_2)=1$.
$$C=\left(\matrix{0 &1 & 0\cr 0  &0 &1\cr -1&0 &0\cr}\right)\Longrightarrow
{\rm (H-4)'}\;\left(\matrix{ 0 &0 & 0\cr -e_3 &e_2 &0\cr 0& 0
&0\cr}\right)\cong {\rm (H-4)}\;\;\left(\matrix{ e_1 &e_3 & 0\cr 0
&0 &0\cr 0& 0 &0\cr}\right).$$ \indent (BII-3) $f_{11}(e_1)=
f_{33}(e_1)=f_{11}(e_2)=f_{33}(e_2)=0$.
$$C=\left(\matrix{0 &1 & 0\cr 0 &0 &1\cr -1&0 &0\cr}\right)\Longrightarrow
{\rm (H-5)}\;\;\left(\matrix{ 0 &0 & 0\cr -e_3 &0 &0\cr 0& 0
&0\cr}\right).$$

\noindent Case (BIII):\quad $f_{11}(e_1)=f_{11}(e_2)=0$:
$$C=\left(\matrix{0 &1 & 0\cr 0 & 0&1\cr -1&0 &0\cr}\right)\Longrightarrow
{\rm (H-6)}\;\;\left(\matrix{ 0 &0 & 0\cr -e_3&e_1 &0\cr 0& 0
&0\cr}\right).$$

\noindent Case (BIV):\quad $f_{11}(e_1)=f_{11}(e_2)=0$:
$$C_{\lambda}=\left(\matrix{0 &1 & 0\cr 0 & 1&\lambda\cr 1&0
&0\cr}\right),\lambda\ne 0, \Longrightarrow {\rm
(H-7)}_{\lambda}\;\;\left(\matrix{ e_3 &e_3 & 0\cr 0&\lambda e_3
&0\cr 0& 0 &0\cr}\right),\lambda\ne 0.$$
$$C=\left(\matrix{0 &0& -\frac{1}{2}\cr 0 & \frac{1}{2}&0\cr 1&0
&0\cr}\right)\Longrightarrow {\rm (H-8)}\;\;\left(\matrix{
0&\frac{1}{2}e_3 & 0\cr -\frac{1}{2} e_3&0 &0\cr 0& 0
&0\cr}\right).$$

\noindent Case (BV):\quad $f_{11}(e_1)=f_{11}(e_2)=0$:
$$C=\left(\matrix{1 &0 & 0\cr 0 & 0&1\cr 0&1 &0\cr}\right)\Longrightarrow
{\rm (H-9)}\;\;\left(\matrix{ 0 &e_3 & 0\cr 0&e_1 &0\cr 0& 0
&0\cr}\right).$$

\noindent Case (BVI):\quad $f_{11}(e_1)=f_{11}(e_2)=0$:
$$C_{\lambda}=\left(\matrix{0 &1 & 0\cr 0 & 0&\lambda\cr \lambda-1&0
&0\cr}\right),\lambda\ne 0,1\Longrightarrow {\rm
(H-10)_{\lambda}}\;\;\left(\matrix{ 0 &\frac{\lambda}{\lambda-1}e_3
& 0\cr \frac{1}{\lambda-1}e_3&\lambda e_1 &0\cr 0& 0
&0\cr}\right),\lambda\ne 0,1.$$

{\bf Proof}\quad It is as the same as the proof of Proposition 3.2
with the computation case by case. \hfill $\Box$

Moreover, through a direct computation, we know

{\bf Proposition 3.5}\quad With the notations as above, among the
left-symmetric algebras on ${\cal H}$, we have

a) Associative algebras: (H-5), (H-7)$_{\lambda}$($\lambda\ne0$),
(H-8);

b) Transitive left-symmetric algebras:(H-5), (H-6),
(H-7)$_{\lambda}$($\lambda\ne0$), (H-8), (H-9), (H-10)$_\lambda$
($\lambda\ne 0,1$);

c) Novikov algebras: (H-1), (H-2), (H-5), (H-6),
(H-7)$_{\lambda}$($\lambda\ne0$), (H-8), (H-9), (H-10)$_\lambda$
($\lambda\ne 0,1$);

d) Bi-symmetric algebras:(H-5), (H-6),
(H-7)$_{\lambda}$($\lambda\ne0$), (H-8), (H-9), (H-10)$_\lambda$
($\lambda\ne 0,1$);

e) There is not any simple left-symmetric algebra on ${\cal H}$.

\subsection{The left-symmetric algebras on ${\cal N}$}

The automorphism group of ${\cal N}$ is
$${\rm Aut}({\cal N})=\{
\left( \matrix{ a_{11} & 0 & 0\cr  0 & a_{22} & 0\cr a_{31} &
a_{32}& 1\cr}\right)|a_{11}a_{22}\ne 0\}.\eqno (3.5)$$

 {\bf
Proposition 3.6}\quad The equivalent classes of the representations
of ${\cal N}$ are divided into the following cases (maybe there are
equivalent classes in the same case, but not equivalent in different
cases): Set
$$B_J=\left(\matrix{ 0&0&0\cr
0&0&0\cr 1&0&0\cr}\right),\;C_J=\left(\matrix{ 0&0&0\cr
1&0&0\cr 0&1&0\cr}\right)$$
${\rm (AI)}\; f(e_1)=\left(\matrix{f_{11} (e_1) & 0 &0\cr 0
&f_{22}(e_1) &0\cr
0 & 0&f_{33}(e_1)\cr}\right), f(e_2)=0, f(e_3)=\left(\matrix{ f_{11}(e_3)&
0 &0\cr 0 &f_{22}(e_3) &0\cr
0 & 0&f_{33}(e_3)\cr}\right).$

\noindent ${\rm (AII)}\; f(e_1)=\left(\matrix{f_{11} (e_1) & 0 &0\cr 0
&f_{11}(e_1) &0\cr
0 & 0&f_{33}(e_1)\cr}\right), f(e_2)=0, f(e_3)=\left(\matrix{ f_{11}(e_3)&
0 &0\cr 1 &f_{11}(e_3) &0\cr
0 & 0&f_{33}(e_3)\cr}\right).$

\noindent${\rm (AIII)}\; f(e_1)=\left(\matrix{f_{11} (e_1) & 0 &0\cr 0
&f_{11}(e_1) &0\cr
0 & 0&f_{11}(e_1)\cr}\right), f(e_2)=0, f(e_3)=\left(\matrix{ f_{11}(e_3)&
0 &0\cr 1&f_{11}(e_3) &0\cr
0 & 1&f_{11}(e_3)\cr}\right).$

\noindent${\rm (AIV)}\; f(e_1)=\left(\matrix{f_{11} (e_1) & 0 &0\cr 1
&f_{11}(e_1) &0\cr
0 & 0&f_{33}(e_1)\cr}\right), f(e_2)=0, f(e_3)=\left(\matrix{ f_{11}(e_3)&
0 &0\cr 0&f_{11}(e_3) &0\cr
0 & 0&f_{33}(e_3)\cr}\right).$

\noindent${\rm (AV)}\; f(e_1)=\left(\matrix{f_{11} (e_1) & 0 &0\cr 1
&f_{11}(e_1) &0\cr
0 & 0&f_{11}(e_1)\cr}\right), f(e_2)=0, f(e_3)=\left(\matrix{ f_{11}(e_3)&
0 &0\cr 0&f_{11}(e_3) &0\cr
1 & 0&f_{11}(e_3)\cr}\right).$

\noindent${\rm (AVI)}\; f(e_1)=\left(\matrix{f_{11} (e_1) & 0 &0\cr 0
&f_{11}(e_1) &0\cr
0 & 1&f_{11}(e_1)\cr}\right), f(e_2)=0, f(e_3)=\left(\matrix{ f_{11}(e_3)&
0 &0\cr 0&f_{11}(e_3) &0\cr
1 & 0&f_{11}(e_3)\cr}\right).$

\noindent${\rm (AVII)}\; f(e_1)=\left(\matrix{f_{11} (e_1) & 0 &0\cr 0
&f_{11}(e_1) &0\cr
1 & 0&f_{11}(e_1)\cr}\right), f(e_2)=0, f(e_3)=\left(\matrix{ f_{11}(e_3)&
0 &0\cr 1&f_{11}(e_3) &0\cr
0& 1&f_{11}(e_3)\cr}\right).$

\noindent${\rm (AVIII)}\; f(e_1)=\left(\matrix{f_{11} (e_1) & 0 &0\cr 1
&f_{11}(e_1) &0\cr
0 & 1&f_{11}(e_1)\cr}\right), f(e_2)=0, f(e_3)=\left(\matrix{ f_{11}(e_3)&
0 &0\cr 0&f_{11}(e_3) &0\cr
0& 0&f_{11}(e_3)\cr}\right).$

\noindent${\rm (AIX)}\; f(e_1)=\left(\matrix{f_{11} (e_1) & 0 &0\cr 1
&f_{11}(e_1) &0\cr
0 & 1&f_{11}(e_1)\cr}\right), f(e_2)=0, f(e_3)=\left(\matrix{ f_{11}(e_3)&
0 &0\cr 0&f_{11}(e_3) &0\cr
1 & 0&f_{11}(e_3)\cr}\right).$

\noindent${\rm (BI)}\; f(e_1)=\left(\matrix{f_{11} (e_1) & 0 &0\cr 0
&f_{22}(e_1) &0\cr
0 & 0&f_{11}(e_1)\cr}\right), f(e_2)=B_J, f(e_3)=\left(\matrix{
f_{33}(e_3)+1&
0 &0\cr 0 &f_{22}(e_3) &0\cr
0 & 0&f_{33}(e_3)\cr}\right).$

\noindent{\small ${\rm (BII)}\; f(e_1)=\left(\matrix{f_{11} (e_1) & 0 &0\cr
0
&f_{11}(e_1) &0\cr
0 & 0&f_{11}(e_1)\cr}\right), f(e_2)=B_J, f(e_3)=\left(\matrix{
f_{33}(e_3)+1&
0 &0\cr 1 &f_{33}(e_3)+1 &0\cr
0 & 0&f_{33}(e_3)\cr}\right).$

\noindent${\rm (BIII)}\; f(e_1)=\left(\matrix{f_{11} (e_1) & 0 &0\cr 0
&f_{11}(e_1) &0\cr
0 & 0&f_{11}(e_1)\cr}\right), f(e_2)=B_J,
f(e_3)=\left(\matrix{ f_{33}(e_3)+1&
0 &0\cr 0&f_{33}(e_3) &0\cr
0 & 1&f_{33}(e_3)\cr}\right).$

\noindent${\rm (BIV)}\; f(e_1)=\left(\matrix{f_{11} (e_1) & 0 &0\cr 1
&f_{11}(e_1) &0\cr
0 & 0&f_{11}(e_1)\cr}\right), f(e_2)=B_J, f(e_3)=\left(\matrix{
f_{33}(e_3)+1&
0 &0\cr 0&f_{33}(e_3)+1 &0\cr
0 & 0&f_{33}(e_3)\cr}\right).$

\noindent${\rm (BV)}\; f(e_1)=\left(\matrix{f_{11} (e_1) & 0 &0\cr 0
&f_{11}(e_1) &0\cr
0 & 1&f_{11}(e_1)\cr}\right), f(e_2)=B_J, f(e_3)=\left(\matrix{
f_{33}(e_3)+1&
0 &0\cr 0&f_{33}(e_3)&0\cr
0 & 0&f_{33}(e_3)\cr}\right).$

\noindent${\rm (CI)}\; f(e_1)=\left(\matrix{f_{11} (e_1) & 0 &0\cr 0
&f_{11}(e_1) &0\cr
0 & 0&f_{11}(e_1)\cr}\right), f(e_2)=C_J, f(e_3)=\left(\matrix{
f_{33}(e_3)+2&
0 &0\cr 0&f_{33}(e_3)+1&0\cr
0 & 0&f_{33}(e_3)\cr}\right).$}

{\bf Proposition 3.7}\quad For a representation of ${\cal N}$ given
in the above cases respectively, there exist bijective 1-cocycles if
and only if it is equivalent to one of the following corresponding
cases:

\noindent Case (AI):\quad (AI-1)
$f_{11}(e_1)=f_{22}(e_1)=f_{33}(e_1)=0,f_{11}(e_3)=1,
f_{22}(e_3)=0, f_{33}(e_3)=\lambda, \lambda\in {\bf C}$.
$$\forall \lambda\in {\bf C},\;\; C=\left(\matrix{ 0&1&0\cr 1&0&0\cr
0&0&1\cr}\right)\Longrightarrow
{\rm (N-1)}_{\lambda}=\left(\matrix{ 0&0&0\cr 0&0&0\cr 0&e_2&\lambda
e_3\cr}\right),\lambda\in {\bf C}.$$

(AI-2) $f_{11}(e_1)=0, f_{22}(e_1)=1, f_{33}(e_1)=\lambda, f_{11}(e_3)=1,
f_{22}(e_3)=0, f_{33}(e_3)=\mu, \lambda\ne 0, \mu\ne 0$.
$$\forall \lambda\ne 0,\mu\ne 0,\;\; C=\left(\matrix{ 0&1&1\cr 1&0&0\cr
0&0&\frac{\mu}{\lambda}\cr}\right)\Longrightarrow
{\rm (N-2)}_{\lambda,\mu}=\left(\matrix{
e_1+\frac{\lambda(\lambda-1)}{\mu}e_3&0&\lambda e_3\cr 0&0&0\cr \lambda
e_3&e_2&\mu e_3\cr}\right),\lambda\ne0, \mu\ne 0.$$

(AI-3) $f_{11}(e_1)=0, f_{22}(e_1)=1, f_{33}(e_1)=0, f_{11}(e_3)=1,
f_{22}(e_3)=0, f_{33}(e_3)=\mu, \mu\ne 0$.
$$\forall \mu\ne 0,\;\; C=\left(\matrix{ 0&1&0\cr 1&0&0\cr
0&0&1\cr}\right)\Longrightarrow
{\rm (N-3)}_{\mu}=\left(\matrix{ e_1&0&0\cr 0&0&0\cr 0 &e_2&\mu
e_3\cr}\right),\mu\ne 0.$$

(AI-4) $f_{11}(e_1)=0, f_{22}(e_1)=1, f_{33}(e_1)=0, f_{11}(e_3)=1,
f_{22}(e_3)=f_{33}(e_3)=0$.
$$C_\lambda=\left(\matrix{ 0 &1&-\lambda\cr 1&0&0\cr
0&0&1\cr}\right),\lambda\in {\bf C}\Longrightarrow
{\rm (N-4)}_{\lambda}=\left(\matrix{ e_1+\lambda e_3&0&0\cr 0&0&0\cr 0
&e_2&0\cr}\right),\lambda\in {\bf C}.$$

\noindent Case (AII):\quad (AII-1)
$f_{11}(e_1)=f_{33}(e_1)=0,f_{11}(e_3)=0,f_{33}(e_3)=1$.
$$C=\left(\matrix{ 1&0&0\cr 0&0&1\cr 0&1&0\cr}\right)\Longrightarrow
{\rm (N-5)}=\left(\matrix{ 0&0&0\cr 0&0&0\cr 0&e_2&e_1\cr}\right).$$

(AII-2) $f_{11}(e_1)=f_{33}(e_1)=0,f_{11}(e_3)=1,f_{33}(e_3)=0$.
$$C=\left(\matrix{ 0&0&1\cr 1&0&0\cr 0&1&0\cr}\right)\Longrightarrow
{\rm (N-6)}=\left(\matrix{ 0&0&0\cr 0&0&0\cr 0&e_2&e_3+e_2\cr}\right).$$

(AII-3) $f_{11}(e_1)=0, f_{33}(e_1)=1,f_{11}(e_3)=1,f_{33}(e_3)=0$.
$$C=\left(\matrix{ 0&0&1\cr 1&0&0\cr 0&1&0\cr}\right)\Longrightarrow
{\rm (N-7)}=\left(\matrix{ e_1&0&0\cr 0&0&0\cr 0&e_2&e_3+e_2\cr}\right).$$

(AII-4) $f_{11}(e_1)=1, f_{33}(e_1)=0,f_{11}(e_3)=0,f_{33}(e_3)=1$.
$$C=\left(\matrix{ 0&1&0\cr 0&0&1\cr 1&0&0\cr}\right)\Longrightarrow
{\rm (N-8)}=\left(\matrix{ e_1&0&e_3\cr 0&0&0\cr e_3&e_2&0\cr}\right).$$

\noindent Case (AIII): there does not exist any bijective 1-cocycle;

\noindent Case (AIV):\quad (AIV-1)
$f_{11}(e_1)=f_{33}(e_1)=0,f_{11}(e_3)=\lambda,f_{33}(e_3)=1, \lambda\ne 0$.
$$\forall \lambda\ne 0, \;\;C=\left(\matrix{ 1&0&0\cr 0&0&1\cr
0&\lambda&0\cr}\right)\Longrightarrow
{\rm (N-9)}_{\lambda}=\left(\matrix{ 0&0&\lambda e_1\cr 0&0&0\cr \lambda
e_1&e_2&\lambda e_3\cr}\right),\lambda\ne0$$

(AIV-2) $f_{11}(e_1)=f_{33}(e_1)=0,f_{11}(e_3)=0,f_{33}(e_3)=1$.
$$C=\left(\matrix{ 0&1&0\cr 0&0&1\cr 1&0&0\cr}\right)\Longrightarrow
{\rm (N-10)}=\left(\matrix{ e_3&0&0\cr 0&0&0\cr 0&e_2&0\cr}\right).$$

(AIV-3) $f_{11}(e_1)=f_{33}(e_1)=0,f_{11}(e_3)=1,f_{33}(e_3)=0$.
$$C=\left(\matrix{ 1&0&1\cr 1&0&0\cr 0&1&0\cr}\right)\Longrightarrow
{\rm (N-11)}=\left(\matrix{ 0&0&e_2\cr 0&0&0\cr e_2&e_2&e_3\cr}\right).$$
$$C=\left(\matrix{ 1&0&1\cr 1&0&0\cr 0&1&1\cr}\right)\Longrightarrow
{\rm (N-12)'}=\left(\matrix{ 0&0&e_2\cr 0&0&0\cr
e_2&e_2&e_3+e_2-e_1\cr}\right)\cong
{\rm (N-12)}=\left(\matrix{ 0&0&e_2\cr 0&0&0\cr
e_2&e_2&e_3-e_2\cr}\right).$$

(AIV-4) $f_{11}(e_1)=0, f_{33}(e_1)=1,f_{11}(e_3)=1,f_{33}(e_3)=\lambda,
\lambda\in {\bf C}$.
\begin{eqnarray*}
\forall \lambda\in {\bf C},\;\; C=\left(\matrix{ 1&0&1\cr 1&0&0\cr
0&1&\lambda\cr}\right)\Longrightarrow
{\rm (N-13)'}_\lambda &=& \left(\matrix{ e_1-e_2&0&\lambda
e_1+(1-\lambda)e_2\cr 0&0&0\cr \lambda e_1+(1-\lambda)e_2
&e_2&e_3+(\lambda^2-\lambda)(e_1-e_2)\cr}\right)\\
\cong {\rm (N-13)}_\lambda &=& \left(\matrix{ e_1-e_2&0&e_2\cr 0&0&0\cr
e_2&e_2&e_3-\lambda e_2\cr}\right),\lambda\in {\bf C}
\end{eqnarray*}

(AIV-5) $f_{11}(e_1)=1, f_{33}(e_1)=0,f_{11}(e_3)=\lambda,f_{33}(e_3)=1,
\lambda\ne0$.
\begin{eqnarray*}
\forall \lambda\ne 0,\;\; C=\left(\matrix{ 0&1&0\cr 0&0&1\cr -\lambda
&\lambda &0\cr}\right)\Longrightarrow
{\rm (N-14)'}_\lambda &=& \left(\matrix{ 2e_1-\frac{1}{\lambda}e_3&0&\lambda
e_1\cr 0&0&0\cr \lambda e_1
&e_2&\lambda e_3\cr}\right)\\
\cong {\rm (N-14)}_\lambda &=& \left(\matrix{
e_1-\frac{1}{\lambda}e_3&0&e_3\cr 0&0&0\cr e_3&e_2&0\cr}\right),\lambda\ne0
\end{eqnarray*}

\noindent Case (AV): there does not exist any bijective 1-cocycle;

\noindent Case (AVI): $f_{11}(e_1)=0,f_{11}(e_3)=1$.
$$C=\left(\matrix{ 0&1&0\cr 1&0&0\cr 0&0&1\cr}\right)\Longrightarrow
{\rm (N-15)}=\left(\matrix{ 0&0&e_1\cr 0&0&0\cr
e_1&e_2&e_3+e_2\cr}\right).$$

\noindent Case (AVII): there does not exist any bijective 1-cocycle;

\noindent Case (AVIII): $f_{11}(e_1)=0,f_{11}(e_3)=1$.
$$C=\left(\matrix{ 0&1&0\cr 1&0&0\cr 0&0&1\cr}\right)\Longrightarrow
{\rm (N-16)}=\left(\matrix{ e_2&0&e_1\cr 0&0&0\cr e_1&e_2&e_3\cr}\right).$$

\noindent Case (AIX): $f_{11}(e_1)=0,f_{11}(e_3)=1$.
$$C=\left(\matrix{ 0&1&0\cr 1&0&0\cr 0&0&1\cr}\right)\Longrightarrow
{\rm (N-17)}=\left(\matrix{ e_2&0&e_1\cr 0&0&0\cr
e_1&e_2&e_3+e_2\cr}\right).$$

\noindent Case (BI):\quad (BI-1) $f_{11}(e_1)=f_{22}(e_1)=0,f_{22}(e_3)=0,
f_{33}(e_3)=\lambda,\lambda\ne 0$.
$$\forall \lambda\ne 0,\;\; C=\left(\matrix{ 0&1&0\cr 1&0&0\cr
0&0&\lambda\cr}\right)\Longrightarrow
{\rm (N-18)}_{\lambda}=\left(\matrix{ 0&0&0\cr 0&0&\lambda e_2\cr
0&(\lambda+1) e_2&\lambda e_3\cr}\right),\lambda\ne 0;$$
when $\lambda =-1$, there is an additional equivalent class:
$$C=\left(\matrix{ 0&1&0\cr 1&0&0\cr 1&0&-1\cr}\right)\Longrightarrow
{\rm (N-19)}=\left(\matrix{ 0&0&0\cr 0&0&-e_2\cr 0&0&-e_3+e_2\cr}\right);$$
when $\lambda =1$, there is an additional equivalent class:
$$C=\left(\matrix{ 0&1&0\cr 0&0&1\cr 1&0&0\cr}\right)\Longrightarrow
{\rm (N-20)}=\left(\matrix{ 0&0&0\cr 0&e_3&0\cr 0&e_2&2e_3\cr}\right).$$

(BI-2) $f_{11}(e_1)=f_{22}(e_1)=0,f_{22}(e_3)=1,
f_{33}(e_3)=-1$.
$$C=\left(\matrix{ 1&0&0\cr -1&1&0\cr 0&0&1\cr}\right)\Longrightarrow
{\rm (N-21)}=\left(\matrix{ 0&0&0\cr 0&0&e_1\cr 0&e_2+e_1&-e_3\cr}\right);$$

(BI-3) $f_{11}(e_1)=1, f_{22}(e_1)=\lambda,f_{22}(e_3)=\mu,
f_{33}(e_3)=0,\lambda\ne 0,\mu\ne0$.
$$\forall \lambda\ne0, \mu\ne 0,\;\; C=\left(\matrix{ 0&1&1\cr 1&0&0\cr
0&\frac{\mu}{\lambda}&0\cr}\right)\Longrightarrow
{\rm (N-22)}_{\lambda,\mu}=\left(\matrix{
e_1+\frac{\lambda(\lambda-1)}{\mu}e_3&e_2&\lambda e_3\cr e_2&0&0\cr
\lambda e_3&e_2&\mu e_3\cr}\right),\lambda\ne 0, \mu\ne 0.$$

(BI-4) $f_{11}(e_1)=1, f_{22}(e_1)=0, f_{22}(e_3)=\mu,
f_{22}(e_3)=0, \mu\ne 0,1$.
$$\forall \mu\ne 0,1,\;\; C=\left(\matrix{ 0&0&1\cr 1&0&0\cr
0&1&0\cr}\right)\Longrightarrow
{\rm (N-23)}_{\mu}=\left(\matrix{ e_1&e_2&0\cr e_2&0&0\cr 0 &e_2&\mu
e_3\cr}\right),\mu\ne 0,1.$$

(BI-5) $f_{11}(e_1)=1, f_{22}(e_1)=0, f_{22}(e_3)=f_{33}(e_3)=0$.
$$C_{\lambda}=\left(\matrix{ 0&-\lambda&1\cr 1&0&0\cr
0&1&0\cr}\right),\;\lambda\in {\bf C}\Longrightarrow
{\rm (N-24)}_{\lambda}=\left(\matrix{ e_1+\lambda e_3&e_2&0\cr e_2&0&0\cr 0
&e_2&0\cr}\right),\lambda\in {\bf C}.$$

(BI-6) $f_{11}(e_1)=1, f_{22}(e_1)=0,  f_{22}(e_3)=1, f_{33}(e_3)=0$.
$$C_{\lambda}=\left(\matrix{ 0&0&1\cr 1&-\lambda&0\cr
0&1&0\cr}\right),\;\lambda\in {\bf C}\Longrightarrow
{\rm (N-25)}_{\lambda}=\left(\matrix{ e_1&e_2+\lambda e_3&0\cr e_2+\lambda
e_3&0&0\cr 0 &e_2&e_3\cr}\right),\lambda\in {\bf C}.$$

(BI-7) $f_{11}(e_1)=0, f_{22}(e_1)=1, f_{22}(e_3)=0,
f_{33}(e_3)=\lambda, \lambda\ne 0$.
$$\forall \lambda\ne 0\;\; C=\left(\matrix{ 0&1&0\cr 1&0&0\cr
0&0&\lambda\cr}\right)\Longrightarrow
{\rm (N-26)}_{\lambda}=\left(\matrix{ e_1&0&0\cr 0&0&\lambda e_2\cr 0
&(\lambda+1)e_2&\lambda e_3\cr}\right),\lambda\ne 0;$$
when $\lambda =-1$, there is three additional equivalent classes:
$$C=\left(\matrix{ 0&1&0\cr 1&0&0\cr 1&0&-1\cr}\right)\Longrightarrow
{\rm (N-27)}=\left(\matrix{ e_1&0&0\cr 0&0&-e_2\cr
0&0&-e_3+e_2\cr}\right);$$
$$C=\left(\matrix{ -1&1&0\cr 1&0&0\cr 0&0&-1\cr}\right)\Longrightarrow
{\rm (N-28)}=\left(\matrix{ e_1+e_2&0&0\cr 0&0&-e_2\cr
0&0&-e_3\cr}\right);$$
$$C=\left(\matrix{ -1&1&0\cr 1&0&0\cr 1&0&-1\cr}\right)\Longrightarrow
{\rm (N-29)}=\left(\matrix{ e_1+e_2&0&0\cr 0&0&-e_2\cr
0&0&-e_3+e_2\cr}\right);$$
when $\lambda=1$, there is an additional equivalent class:
$$C=\left(\matrix{ 0&1&0\cr 0&0&1\cr 1&0&0\cr}\right)\Longrightarrow
{\rm (N-30)}=\left(\matrix{ e_1&0&0\cr 0&e_3&0\cr 0&e_2&2e_3\cr}\right).$$

\noindent Case (BII): $f_{11}(e_1)=1,f_{33}(e_3)=-1$.
$$C=\left(\matrix{ 0&1&1\cr 1&0&0\cr 1&1&0\cr}\right)\Longrightarrow
{\rm (N-31)}=\left(\matrix{ e_1&e_2&e_3\cr e_2&0&0\cr
e_3&e_2&e_3+e_2\cr}\right).$$

\noindent Case (BIII): $f_{11}(e_1)=1,f_{11}(e_3)=0$.
$$C=\left(\matrix{ 0&0&1\cr 1&0&0\cr 0&1&0\cr}\right)\Longrightarrow
{\rm (N-32)}=\left(\matrix{ e_1&e_2&e_3\cr e_2&0&0\cr
e_3&e_2&0\cr}\right).$$

\noindent Case (BIV): (BIV-1) $f_{11}(e_1)=0, f_{33}(e_3)=-1$.
$$C=\left(\matrix{ 0&1&0\cr 1&0&0\cr 0&0&-1\cr}\right)\Longrightarrow
{\rm (N-33)}=\left(\matrix{ e_2&0&0\cr 0&0&-e_2\cr 0&0&-e_3\cr}\right);$$
$$C=\left(\matrix{ 0&1&0\cr 1&0&0\cr 1&0&-1\cr}\right)\Longrightarrow
{\rm (N-34)}=\left(\matrix{ e_2&0&0\cr 0&0&-e_2\cr
0&0&-e_3+e_2\cr}\right).$$

(BIV-2) $f_{11}(e_1)=0, f_{33}(e_3)=0$.
$$C=\left(\matrix{ 0&0&1\cr 0&1&0\cr 1&0&0\cr}\right)\Longrightarrow
{\rm (N-35)}=\left(\matrix{ 0&e_3&0\cr e_3&0&0\cr 0&e_2&e_3\cr}\right).$$

(BIV-3) $f_{11}(e_1)=0, f_{33}(e_3)=1$.
$$C=\left(\matrix{ 1&0&0\cr 0&0&1\cr 0&2&0\cr}\right)\Longrightarrow
{\rm (N-36)}=\left(\matrix{ 0&0&2e_1\cr 0&e_1&0\cr
2e_1&e_2&2e_3\cr}\right).$$

(BIV-4) $f_{11}(e_1)=1,f_{11}(e_3)=\lambda,\lambda\in {\bf C}$.
\begin{eqnarray*}
&&\forall \lambda\in {\bf C},\;\; C=\left(\matrix{ 0&1&1\cr 1&0&0\cr
-(\lambda+1)&\lambda+1&\lambda\cr}\right)\Longrightarrow\\
&&{\rm (N-37)'}_\lambda = \left(\matrix{ e_1+e_2&e_2& e_3+(1+\lambda)e_2\cr
e_2&0&\lambda e_2\cr e_3+(1+\lambda)e_2
&(\lambda+1)e_2&(2\lambda+1) e_3-\lambda(\lambda+1) (e_1- e_2)\cr}\right)\\
&&\cong {\rm (N-37)}_\lambda = \left(\matrix{ e_1+e_2&e_2&e_2+e_3\cr
e_2&0&0\cr e_2+e_3&e_2&e_3-\lambda e_2\cr}\right),\lambda\in {\bf C}
\end{eqnarray*}

\noindent Case (BV):\quad (BV-1)
$f_{11}(e_1)=0,f_{33}(e_3)=\lambda,\lambda\ne 0$.
$$\forall \lambda\ne 0,\;\; C=\left(\matrix{ 0&1&0\cr 1&0&0\cr
0&0&\lambda\cr}\right)\Longrightarrow
{\rm (N-38)}_{\lambda}=\left(\matrix{ 0&0&\lambda e_1\cr 0&0&\lambda e_2\cr
\lambda e_1&(\lambda+1) e_2&\lambda e_3\cr}\right),\lambda\ne 0;$$
when $\lambda =-1$, there is three additional equivalent classes:
$$C=\left(\matrix{ 0&1&0\cr 1&0&0\cr -1&0&-1\cr}\right)\Longrightarrow
{\rm (N-39)}=\left(\matrix{ 0&0&-e_1\cr 0&0&-e_2\cr
-e_1&0&-e_3+e_2\cr}\right);$$
$$C=\left(\matrix{ 1&1&0\cr 1&0&0\cr 0&0&-1\cr}\right)\Longrightarrow
{\rm (N-40)}=\left(\matrix{ 0&0&-e_1+e_2\cr 0&0&-e_2\cr
-e_1+e_2&0&-e_3\cr}\right);$$
$$C=\left(\matrix{ 1&1&0\cr 1&0&0\cr -1&0&-1\cr}\right)\Longrightarrow
{\rm (N-41)}=\left(\matrix{ 0&0&-e_1+e_2\cr 0&0&-e_2\cr
-e_1+e_2&0&-e_3+e_2\cr}\right);$$
when $\lambda =1$, there is an additional equivalent class:
$$C=\left(\matrix{ 0&1&0\cr 1&1&0\cr 0&0&1\cr}\right)\Longrightarrow
{\rm (N-42)}=\left(\matrix{ 0&0&e_1\cr 0&0&e_2-e_1\cr
e_1&2e_2-e_1&e_3\cr}\right).$$

(BV-2) $f_{11}(e_1)=1,f_{33}(e_3)=\lambda,\lambda\ne 0$.
$$\forall \lambda\ne 0,\;\; C=\left(\matrix{ 0&0&1\cr 1&0&0\cr
0&-\lambda&\lambda\cr}\right)\Longrightarrow
{\rm (N-43)}_\lambda=\left(\matrix{ e_1-\frac{1}{\lambda}e_3&e_2&e_3\cr
e_2&0&0\cr e_3&e_2&0\cr}\right),\lambda\ne 0.$$

\noindent Case (CI): (CI-1) $f_{11}(e_1)=0,f_{33}(e_3)=-2$.
$$C=\left(\matrix{ 1&0&0\cr 0&1&0\cr 0&0&-2\cr}\right)\Longrightarrow
{\rm (N-44)}=\left(\matrix{ 0&0&0\cr 0&e_1&-2e_2\cr
0&-e_2&-2e_3\cr}\right).$$

(CI-2) $f_{11}(e_1)=1,f_{33}(e_3)=0$.
$$C=\left(\matrix{ 1&0&1\cr 0&1&0\cr 1&0&0\cr}\right)\Longrightarrow
{\rm (N-45)}=\left(\matrix{ e_1&e_2&e_3\cr e_2&e_3&0\cr
e_3&e_2&2e_3\cr}\right).$$

{\bf Remark}\quad It is easy to see that we can extend the extent of
some parameters appearing in above left-symmetric algebras:

(1) $ {\rm (N-2)}_{\lambda=0,\mu\ne 0}\cong {\rm
(N-3)}_{\mu\ne0}$;\quad ${\rm (N-3)}_{\mu=0}\cong {\rm
(N-4)}_{\lambda=0}\cong {\rm (N-26)}_{\lambda=0}$;

(2) $ {\rm (N-22)}_{\lambda=0,\mu\ne 0,1}\cong {\rm (N-23)}_{\mu\ne
0,1}$;\quad ${\rm (N-22)}_{\lambda=0, \mu=1}\cong {\rm
(N-23)}_{\mu=1}\cong {\rm (N-25)}_{\lambda=0}$;

\mbox{}\hspace{0.6cm} ${\rm (N-23)}_{\mu=0}\cong {\rm
(N-24)}_{\lambda=0}$;

(3) ${\rm (N-1)}_{\mu=0}\cong {\rm (N-9)}_{\lambda=0}\cong {\rm
(N-18)}_{\lambda=0}\cong {\rm (N-38)}_{\lambda=0}$.

{\bf Proposition 3.8}\quad With the notations as above, among the
left-symmetric algebras on ${\cal N}$, we have

a) Associative algebras: (N-1)$_1$, (N-3)$_1$, (N-9)$_1$,  (N-18)$_{-1}$,
(N-22)$_{1,1}$, (N-26)$_{-1}$;

b) Transitive left-symmetric algebras:  (N-1)$_0$, (N-5), (N-10);

c) Novikov algebras: (N-1)$_0$, (N-4)$_0$, (N-5),
(N-18)$_\lambda$($\lambda\ne 0$), (N-19),
(N-26)$_\lambda$($\lambda\ne 0$), (N-27),  (N-38)$_\lambda$($\lambda\ne 0$),
(N-39);

d) Bi-symmetric algebras:(N-1)$_1$, (N-3)$_1$, (N-6), (N-7), (N-9)$_1$,
(N-15), (N-18)$_{-1}$, (N-19), (N-22)$_{1,1}$,
(N-26)$_{-1}$; (N-27), (N-31), (N-38)$_{-1}$, (N-39).

e) There is not any simple left-symmetric algebra on ${\cal N}$. But (N-30)
is semisimple.

\subsection{The left-symmetric algebras on ${\cal D}_1$}

The automorphism group of ${\cal D}_1$ is
$${\rm Aut}({\cal D}_1)=\{
\left( \matrix{ a_{11} & a_{12} & 0\cr  a_{21} & a_{22} & 0\cr
a_{31} & a_{32}& 1\cr}\right)|a_{11}a_{22}-a_{12}a_{21}\ne 0\}.\eqno
(3.6)$$

{\bf Proposition 3.9}\quad The equivalent classes of the
representations of ${\cal D}_1$ are divided into the following
cases:

\noindent ${\rm (AI)}\; f(e_1)=0, f(e_2)=0, f(e_3)=\left(\matrix{
f_{11}(e_3)&
0 &0\cr 0 &f_{22}(e_3) &0\cr
0 & 0&f_{33}(e_3)\cr}\right).$

\noindent ${\rm (AII)}\; f(e_1)=0, f(e_2)=0, f(e_3)=\left(\matrix{
f_{11}(e_3)&
0 &0\cr 1 &f_{11}(e_3) &0\cr
0 & 0&f_{33}(e_3)\cr}\right).$

\noindent ${\rm (AIII)}\; f(e_1)=0, f(e_2)=0, f(e_3)=\left(\matrix{
f_{11}(e_3)&
0 &0\cr 1 &f_{11}(e_3) &0\cr
0 & 1&f_{11}(e_3)\cr}\right).$

\noindent ${\rm (AIV)}\; f(e_1)=0, f(e_2)=\left(\matrix{0 & 0 &0\cr 1
&0 &0\cr
0 & 0&0\cr}\right), f(e_3)=\left(\matrix{ f_{22}(e_3)+1&
0 &0\cr 0&f_{22}(e_3) &0\cr
0 & 0&f_{33}(e_3)\cr}\right).$

\noindent ${\rm (AV)}\; f(e_1)=0, f(e_2)=\left(\matrix{0 & 0 &0\cr 1
&0 &0\cr
0 & 0&0\cr}\right), f(e_3)=\left(\matrix{ f_{33}(e_3)&
0 &0\cr 0&f_{33}(e_3)-1 &0\cr
1 & 0&f_{33}(e_3)\cr}\right).$

\noindent ${\rm (AVI)}\; f(e_1)=0, f(e_2)=\left(\matrix{0 & 0 &0\cr 0
&0 &0\cr
0 & 1&0\cr}\right), f(e_3)=\left(\matrix{ f_{33}(e_3)&
0 &0\cr 0&f_{33}(e_3)+1 &0\cr
1 & 0&f_{33}(e_3)\cr}\right).$

\noindent ${\rm (AVII)}\; f(e_1)=0, f(e_2)=\left(\matrix{0 & 0 &0\cr 1
&0 &0\cr
0 & 1&0\cr}\right), f(e_3)=\left(\matrix{ f_{33}(e_3)+2&
0 &0\cr 0&f_{33}(e_3)+1 &0\cr
0& 0&f_{33}(e_3)\cr}\right).$

\noindent ${\rm (BI)}\; f(e_1)=\left(\matrix{ 0&0&0\cr 1&0&0\cr
0&0&0\cr}\right), f(e_2)=\left(\matrix{0 & 0 &0\cr 0
&0 &0\cr
1 &0 &0\cr}\right), f(e_3)=\left(\matrix{ f_{33}(e_3)+1&
0 &0\cr 0&f_{33}(e_3) &0\cr
0 & 0&f_{33}(e_3)\cr}\right).$

\noindent ${\rm (BII)}\; f(e_1)=\left(\matrix{ 0&0&0\cr 0&0&0\cr
0&1&0\cr}\right), f(e_2)=\left(\matrix{0 & 0 &0\cr 0
&0 &0\cr
1 &0 &0\cr}\right), f(e_3)=\left(\matrix{ f_{33}(e_3)+1&
0 &0\cr 0&f_{33}(e_3)+1 &0\cr
0 & 0&f_{33}(e_3)\cr}\right).$

{\bf Proposition 3.10}\quad  For a representation of ${\cal D}_1$
given in the above cases respectively, there exist bijective
1-cocycles if and only if it is equivalent to one of the following
corresponding cases:

\noindent Case (AI):\quad $f_{11}(e_3)=f_{22}(e_3)=1,
f_{33}(e_3)=\lambda,\lambda\in {\bf C}$.
$$\forall \lambda\in {\bf C},\;\; C=\left(\matrix{ 1&0&0\cr 0&1&0\cr
0&0&1\cr}\right)\Longrightarrow
{{\rm ({\bar D}_1-1)}}_\lambda = \left(\matrix{ 0&0&0\cr 0&0&0\cr  e_1
&e_2&\lambda e_3\cr}\right), \lambda\in {\bf C}.$$

\noindent Case (AII):\quad $f_{11}(e_3)=f_{33}(e_3)=1$.
$$C=\left(\matrix{ 0&0&1\cr 1&0&0\cr 0&1&0\cr}\right)\Longrightarrow
{\rm ({\bar D}_1-2)} = \left(\matrix{ 0&0&0\cr 0&0&0\cr  e_1
&e_2&e_3+e_2\cr}\right).$$

\noindent Case (AIII): there does not exist any bijective 1-cocycle;

\noindent Case (AIV): \quad $f_{22}(e_3)=\lambda, f_{33}(e_3)=1,\lambda\ne
0$.
$$\forall \lambda\ne 0,\;\;C=\left(\matrix{ 0&0&1\cr 1&0&0\cr
0&\lambda&0\cr}\right)\Longrightarrow
{\rm ({\bar D}_1-3)}_\lambda = \left(\matrix{ 0&0&0\cr 0&0&\lambda e_2\cr
e_1
&(\lambda+1)e_2&\lambda e_3\cr}\right),\lambda \ne 0;$$
when $\lambda=-1$, there is an additional equivalent class:
$$C=\left(\matrix{ 0&0&1\cr 1&0&0\cr 1&-1&0\cr}\right)\Longrightarrow
{\rm ({\bar D}_1-4)} = \left(\matrix{ 0&0&0\cr 0&0&-e_2\cr  e_1
&0&-e_3+e_2\cr}\right);$$
when $\lambda=1$, there is an additional equivalent class:
$$C=\left(\matrix{ 0&0&1\cr 0&1&0\cr 1&0&0\cr}\right)\Longrightarrow
{\rm ({\bar D}_1-5)} = \left(\matrix{ 0&0&0\cr 0&e_3&0\cr  e_1
&e_2&2e_3\cr}\right).$$

\noindent Case (AV):\quad $f_{33}(e_3)=1$.
$$C=\left(\matrix{ 1&0&0\cr 0&0&1\cr 0&1&0\cr}\right)\Longrightarrow
{\rm ({\bar D}_1-6)} = \left(\matrix{ 0&0&0\cr 0&0&e_1\cr  e_1
&e_1+e_2&0\cr}\right).$$

\noindent Case (AVI):\quad $f_{33}(e_3)=1$.
$$C=\left(\matrix{ 1&0&0\cr 0&1&0\cr 0&0&1\cr}\right)\Longrightarrow
{\rm ({\bar D}_1-7)} = \left(\matrix{ 0&0&0\cr 0&0&e_2\cr  e_1
&2e_2&e_3+e_1\cr}\right).$$

\noindent Case (AVII):\quad $f_{33}(e_3)=-1$.
$$C=\left(\matrix{ 1&0&0\cr 0&1&0\cr 0&0&-1\cr}\right)\Longrightarrow
{\rm ({\bar D}_1-8)} = \left(\matrix{ 0&0&0\cr 0&e_1&-e_2\cr  e_1
&0&-e_3\cr}\right).$$

\noindent Case (BI):\quad $f_{33}(e_3)=1$.
$$C=\left(\matrix{ 1&0&0\cr 0&0&1\cr 0&1&0\cr}\right)\Longrightarrow
{\rm ({\bar D}_1-9)} = \left(\matrix{ 0&0&e_1\cr 0&e_1&0\cr  2e_1
&e_2&e_3\cr}\right);$$
$$C=\left(\matrix{ 0&1&0\cr 0&0&1\cr 1&0&0\cr}\right)\Longrightarrow
{\rm ({\bar D}_1-10)} = \left(\matrix{ e_3&0&0\cr 0&e_3&0\cr  e_1
&e_2&2e_3\cr}\right).$$

\noindent Case (BII):\quad $f_{33}(e_3)=\lambda,\lambda\ne 0$.
$$\forall \lambda\ne 0,\;\;C=\left(\matrix{ 0&1&0\cr 1&0&0\cr
0&0&\lambda\cr}\right)\Longrightarrow
{\rm ({\bar D}_1-11)}_\lambda = \left(\matrix{ 0&0&\lambda e_1\cr
0&0&\lambda e_2\cr  (\lambda+1)e_1
&(\lambda+1)e_2&\lambda e_3\cr}\right),\lambda \ne 0;$$
when $\lambda=-1$, there is an additional equivalent class:
$$C=\left(\matrix{ 0&1&0\cr 1&0&0\cr 1&0&-1\cr}\right)\Longrightarrow
{\rm ({\bar D}_1-12)} = \left(\matrix{ 0&0&-e_1\cr 0&0&-e_2\cr  0
&0&-e_3+e_2\cr}\right).$$

{\bf Remark}\quad It is easy to see that we can extend the extent of
the some parameters appearing in above left-symmetric algebras:
$$ {\rm ({\bar D}_1-1)}_{\lambda=0}\cong {\rm ({\bar
D}_1-3)}_{\lambda=0}\cong {\rm ({\bar D}_1-11)}_{\lambda=0}.$$

{\bf Proposition 3.11}\quad With the notations as above, among the
left-symmetric algebras on ${\cal D}_1$, we have

a) Associative algebras: $({\bar {\rm D}}_1$-1)$_1$, $({\bar {\rm
D}}_1$-11)$_{-1}$;

b) Transitive left-symmetric algebras:$({\bar {\rm D}}_1$-1)$_0$,$({\bar
{\rm D}}_1$-6);

c) Novikov algebras: $({\bar {\rm D}}_1$-1)$_0$, $({\bar {\rm
D}}_1$-11)$_\lambda$ ($\lambda\ne 0$), $({\bar {\rm D}}_1$-12);

d) Bi-symmetric algebras:$({\bar {\rm D}}_1$-1)$_1$, $({\bar {\rm D}}_1$-2),
$({\bar {\rm D}}_1$-11)$_{-1}$, $({\bar {\rm D}}_1$-12);

e) There is one simple left-symmetric algebra on ${\cal D}_1$:$({\bar {\rm D}}_1$-10).

\subsection{ The left-symmetric algebras on ${\cal D}_l, \;0<|l|<1$
or $l=e^{i\theta}, 0<\theta\leq \pi$}

Throughout this subsection, without special saying, $0<|l|<1$ or
$l=e^{i\theta}, 0<\theta\leq \pi$. The automorphism group of ${\cal
D}_l$ is
$$\forall\; 0<|l|<1,\;\;{\rm or}\;\; l=e^{i\theta}, 0<\theta< \pi,\;\;{\rm Aut}({\cal D}_l)=\{
\left( \matrix{ a_{11} & 0 & 0\cr  0 & a_{22} & 0\cr a_{31} &
a_{32}& 1\cr}\right)|a_{11}a_{22}\ne 0\}.\eqno (3.7)$$
$${\rm Aut}({\cal D}_{-1})=\{
\left( \matrix{ a_{11} & 0 & 0\cr  0 & a_{22} & 0\cr a_{31} &
a_{32}& 1\cr}\right) |a_{11}a_{22}\ne 0\}\bigcup \{ \left( \matrix{
0& a_{12} & 0\cr  a_{21}&0 & 0\cr a_{31} & a_{32}&
-1\cr}\right)|a_{12}a_{21}\ne 0\}.\eqno (3.8)$$

{\bf Proposition 3.12}\quad The equivalent classes of the
representations of ${\cal D}_l$ are divided into the following
cases:

\noindent ${\rm (AI)}\; f(e_1)=0, f(e_2)=0, f(e_3)=\left(\matrix{
f_{11}(e_3)&
0 &0\cr 0 &f_{22}(e_3) &0\cr
0 & 0&f_{33}(e_3)\cr}\right).$

\noindent ${\rm (AII)}\; f(e_1)=0, f(e_2)=0, f(e_3)=\left(\matrix{
f_{11}(e_3)&
0 &0\cr 1 &f_{11}(e_3) &0\cr
0 & 0&f_{33}(e_3)\cr}\right).$

\noindent ${\rm (AIII)}\; f(e_1)=0, f(e_2)=0, f(e_3)=\left(\matrix{
f_{11}(e_3)&
0 &0\cr 1 &f_{11}(e_3) &0\cr
0 & 1&f_{11}(e_3)\cr}\right).$

\noindent ${\rm (AIV)}\; f(e_1)=0, f(e_2)=\left(\matrix{0 & 0 &0\cr 1
&0 &0\cr
0 & 0&0\cr}\right), f(e_3)=\left(\matrix{ f_{22}(e_3)+l&
0 &0\cr 0&f_{22}(e_3) &0\cr
0 & 0&f_{33}(e_3)\cr}\right).$

\noindent ${\rm (AV)}\; f(e_1)=0, f(e_2)=\left(\matrix{0 & 0 &0\cr 1
&0 &0\cr
0 & 0&0\cr}\right), f(e_3)=\left(\matrix{ f_{33}(e_3)&
0 &0\cr 0&f_{33}(e_3)-l &0\cr
1 & 0&f_{33}(e_3)\cr}\right).$

\noindent ${\rm (AVI)}\; f(e_1)=0, f(e_2)=\left(\matrix{0 & 0 &0\cr 0
&0 &0\cr
0 & 1&0\cr}\right), f(e_3)=\left(\matrix{ f_{33}(e_3)&
0 &0\cr 0&f_{33}(e_3)+l &0\cr
1 & 0&f_{33}(e_3)\cr}\right).$

\noindent ${\rm (AVII)}\; f(e_1)=0, f(e_2)=\left(\matrix{0 & 0 &0\cr 1
&0 &0\cr
0 & 1&0\cr}\right), f(e_3)=\left(\matrix{ f_{33}(e_3)+2l&
0 &0\cr 0&f_{33}(e_3)+l &0\cr
0& 0&f_{33}(e_3)\cr}\right).$

\noindent ${\rm (BI)}\; f(e_1)=\left(\matrix{ 0&0&0\cr 1&0&0\cr
0&0&0\cr}\right), f(e_2)=\left(\matrix{0 & 0 &0\cr 0
&0 &0\cr
1 &0 &0\cr}\right), f(e_3)=\left(\matrix{ f_{33}(e_3)+l&
0 &0\cr 0&f_{33}(e_3)+l-1 &0\cr
0 & 0&f_{33}(e_3)\cr}\right).$

\noindent ${\rm (BII)}\; f(e_1)=\left(\matrix{ 0&0&0\cr 0&0&0\cr
0&1&0\cr}\right), f(e_2)=\left(\matrix{0 & 0 &0\cr 0
&0 &0\cr
1 &0 &0\cr}\right), f(e_3)=\left(\matrix{ f_{33}(e_3)+l&
0 &0\cr 0&f_{33}(e_3)+1 &0\cr
0 & 0&f_{33}(e_3)\cr}\right).$

\noindent ${\rm (BIII)}_{\frac{1}{2}}\;$ only for $l=\frac{1}{2}$,
$$f(e_1)=\left(\matrix{ 0&0&0\cr 0&0&0\cr 1&0&0\cr}\right),
f(e_2)=\left(\matrix{0 & 0 &0\cr 1
&0 &0\cr
0 &1 &0\cr}\right), f(e_3)=\left(\matrix{ f_{33}(e_3)+1&
0 &0\cr 0&f_{33}(e_3)+\frac{1}{2} &0\cr
0 & 0&f_{33}(e_3)\cr}\right).$$

\noindent ${\rm (BIV)}\;(l\ne -1)\; f(e_1)=\left(\matrix{ 0&0&0\cr 1&0&0\cr
0&0&0\cr}\right), f(e_2)=0, f(e_3)=\left(\matrix{ f_{22}(e_3)+1&
0 &0\cr 0&f_{22}(e_3) &0\cr
0 & 0&f_{33}(e_3)\cr}\right).$

\noindent ${\rm (BV)}\;(l\ne -1)\; f(e_1)=\left(\matrix{ 0&0&0\cr 1&0&0\cr
0&0&0\cr}\right), f(e_2)=0, f(e_3)=\left(\matrix{ f_{33}(e_3)&
0 &0\cr 0&f_{33}(e_3)-1 &0\cr
1 & 0&f_{33}(e_3)\cr}\right).$

\noindent ${\rm (BVI)}\;(l\ne -1)\; f(e_1)=\left(\matrix{ 0&0&0\cr 0&0&0\cr
0&1&0\cr}\right), f(e_2)=0, f(e_3)=\left(\matrix{ f_{33}(e_3)&
0 &0\cr 0&f_{33}(e_3)+1 &0\cr
1 & 0&f_{33}(e_3)\cr}\right).$

\noindent ${\rm (CI)}\; (l\ne -1)\;
f(e_1)=\left(\matrix{ 0&0&0\cr 1&0&0\cr 0&1&0\cr}\right), f(e_2)=0,
f(e_3)=\left(\matrix{ f_{33}(e_3)+2&
0 &0\cr 0&f_{33}(e_3)+1 &0\cr
0 & 0&f_{33}(e_3)\cr}\right).$

{\bf Remark}\quad If for the latter four cases, we extend $l$ to
$l=-1$, then it is easy to see that the following cases of
representations of ${\cal D}_{-1}$ are equivalent:
$${\rm (AIV)}\sim {\rm (BIV)},\;\;{\rm (AV)}\sim {\rm (BV)},\;\;{\rm
(AVI)}\sim {\rm (BVI)},\;\;{\rm (AVII)}\sim {\rm (CI)},$$ by the
following linear isomorphism which is in ${\rm Aut}({\cal D}_{-1})$
$$e_1\rightarrow e_2,\;\; e_2\rightarrow e_1,\;\;e_3\rightarrow -e_3.$$

{\bf Proposition 3.13}\quad For a representation of ${\cal D}_l$
given in the above cases respectively, there exist bijective
1-cocycles if and only if it is equivalent to one of the following
corresponding cases:

\noindent Case (AI):\quad $f_{11}(e_3)=1, f_{22}(e_3)=l,
f_{33}(e_3)=\lambda,\lambda\in {\bf C}$.
$$\forall \lambda\in {\bf C},\;\; C=\left(\matrix{ 1&0&0\cr 0&1&0\cr
0&0&1\cr}\right)\Longrightarrow
{\rm (D}_l-1)_\lambda = \left(\matrix{ 0&0&0\cr 0&0&0\cr  e_1
&le_2&\lambda e_3\cr}\right), \lambda\in {\bf C}.$$

\noindent Case (AII):\quad (AII-1)\quad $f_{11}(e_3)=l, f_{33}(e_3)=1$.
$$C=\left(\matrix{ 0&0&1\cr 1&0&0\cr 0&1&0\cr}\right)\Longrightarrow
{\rm (D}_l-2) = \left(\matrix{ 0&0&0\cr 0&0&0\cr  e_1
&le_2&le_3+e_2\cr}\right).$$

Case(AII-2)\quad $f_{11}(e_3)=1, f_{33}(e_3)=l$.
$$C=\left(\matrix{ 1&0&0\cr 0&0&1\cr 0&1&0\cr}\right)\Longrightarrow
{\rm (D}_l-3) = \left(\matrix{ 0&0&0\cr 0&0&0\cr  e_1
&le_2&e_3+e_1\cr}\right).$$

\noindent Case (AIII): there does not exist any bijective 1-cocycle;

\noindent Case (AIV):\quad (AIV-1)\quad $f_{22}(e_3)=\lambda,
f_{33}(e_3)=1,\lambda\ne 0$.
$$\forall \lambda\ne 0,\;\;C=\left(\matrix{ 0&0&1\cr 1&0&0\cr
0&\lambda&0\cr}\right)\Longrightarrow
{\rm (D}_l-4)_\lambda = \left(\matrix{ 0&0&0\cr 0&0&\lambda e_2\cr  e_1
&(\lambda+l)e_2&\lambda e_3\cr}\right),\lambda \ne 0;$$
when $\lambda=-l$, there is an additional equivalent class:
$$C=\left(\matrix{ 0&0&1\cr l&0&0\cr 1&-l&0\cr}\right)\Longrightarrow
{\rm (D}_l-5) = \left(\matrix{ 0&0&0\cr 0&0&-le_2\cr  e_1
&0&-le_3+e_2\cr}\right);$$
when $\lambda=l$, there is an additional equivalent class:
$$C=\left(\matrix{ 0&0&1\cr 0&1&0\cr 1&0&0\cr}\right)\Longrightarrow
{\rm (D}_l-6) = \left(\matrix{ 0&0&0\cr 0&e_3&0\cr  e_1
&le_2&2le_3\cr}\right).$$

(AIV-2):\quad $f_{22}(e_3)=1-l, f_{33}(e_3)=l$.
$$C=\left(\matrix{ 1&0&0\cr \frac{1}{1-l}&0&1\cr
0&1&0\cr}\right)\Longrightarrow
{\rm (D}_l-7) = \left(\matrix{ 0&0&0\cr 0&0&e_1\cr  e_1
&le_2+e_1&(1-l)e_3\cr}\right).$$

(AIV-3)$_{\frac{1}{2}}$:\quad
(only for $l=\frac{1}{2}$)\quad $f_{22}(e_3)=\frac{1}{2},
f_{33}(e_3)=\lambda,\lambda\in {\bf C}$.
$$\forall \lambda\in {\bf C},\;\;C=\left(\matrix{ 1&0&0\cr 0&1&0\cr
0&0&1\cr}\right)\Longrightarrow
{\rm (D_{\frac{1}{2}}-S-1)}_\lambda = \left(\matrix{ 0&0&0\cr 0&e_1&0\cr
e_1
&\frac{1}{2}e_2&\lambda e_3\cr}\right),\lambda \in {\bf C}.$$
(Notice that when $\lambda=\frac{1}{2}$, it need to add $({\rm
D_{\frac{1}{2}}-7})$;
when $\lambda=1$, it need to add $({\rm D_{\frac{1}{2}}-4})_{\frac{1}{2}}$
and $({\rm D_{\frac{1}{2}}-6})$.

\noindent Case (AV):\quad for $l\ne \frac{1}{2}$, there does not exist any
bijective 1-cocycle; for $l=\frac{1}{2}$:\quad
$f_{33}(e_3)=1$.
$$C=\left(\matrix{ 1&0&0\cr 0&1&0\cr 0&0&1\cr}\right)\Longrightarrow
{\rm (D_{\frac{1}{2}}-S-2)}= \left(\matrix{ 0&0&0\cr 0&e_1&0\cr
e_1 &\frac{1}{2}e_2&e_3+e_1\cr}\right).$$

\noindent Case (AVI):\quad (AVI-1) \quad $f_{33}(e_3)=1$.
$$C=\left(\matrix{ 1&0&0\cr 0&1&0\cr 0&0&1\cr}\right)\Longrightarrow
{\rm (D}_l-8) = \left(\matrix{ 0&0&0\cr 0&0&e_2\cr  e_1
&(1+l)e_2&e_3+e_1\cr}\right);$$
for $l=-1$, there is an additional equivalent class:
$$C=\left(\matrix{ 1&0&0\cr 0&1&0\cr 0&1&1\cr}\right)\Longrightarrow
{\rm (D_{-1}-T-1)} = \left(\matrix{ 0&0&0\cr 0&0&e_2\cr  e_1
&0&e_3+e_1+e_2\cr}\right).$$

(AVI-2)$_{\frac{1}{2}}$:\quad
(only for $l=\frac{1}{2}$)\quad $f_{33}(e_3)=\frac{1}{2}$.
$$C=\left(\matrix{ 0&1&0\cr 1&2&0\cr 0&4&1\cr}\right)\Longrightarrow
{\rm (D_{\frac{1}{2}}-S-3)}= \left(\matrix{ 0&0&0\cr 0&e_1&0\cr  e_1
&\frac{1}{2}e_2+e_1&\frac{1}{2}e_3+e_2\cr}\right).$$

\noindent Case (AVII):\quad for $l=\frac{1}{2}$, there does not exist any
bijective 1-cocycle; for
$l\ne\frac{1}{2}$:\quad $f_{33}(e_3)=1-2l$.
$$C=\left(\matrix{ 1&0&0\cr 0&1&0\cr 0&0&1-2l\cr}\right)\Longrightarrow
{\rm (D}_{l\ne {\frac{1}{2}}}-{\rm S}-1)= \left(\matrix{ 0&0&0\cr 0&e_1&
(1-2l)e_2\cr
e_1&(1-l)e_2&(1-2l)e_3\cr}\right);$$
for $l=\frac{1}{3}$, there is an additional equivalent class:
$$C=\left(\matrix{ 1&0&0\cr \frac{3}{2}&0&1\cr
0&1&0\cr}\right)\Longrightarrow
{\rm (D_{\frac{1}{3}}-T-1)} = \left(\matrix{ 0&0&0\cr 0&e_3&e_1\cr  e_1
&\frac{1}{3}e_2+e_1&\frac{2}{3}e_3\cr}\right).$$

\noindent Case (BI):\quad (BI-1)\quad $f_{33}(e_3)=2-l$.
$$C=\left(\matrix{ 0&1&0\cr 1&0&0\cr 0&0&2-l\cr}\right)\Longrightarrow
{\rm (D}_l-9) = \left(\matrix{ e_2&0&0\cr 0&0&(2-l)e_2\cr  e_1
&2e_2&(2-l)e_3\cr}\right).$$

(BI-2)\quad $f_{33}(e_3)=1$.
$$C=\left(\matrix{ 0&0&1\cr 0&1&0\cr 1&0&0\cr}\right)\Longrightarrow
{\rm (D}_l-10) = \left(\matrix{ 0&e_3&0\cr e_3&0&0\cr  e_1
&le_2&(l+1)e_3\cr}\right).$$

(BI-3)$_{l\ne \frac{1}{2}}$:\quad
(only for $l\ne \frac{1}{2}$)\quad $f_{33}(e_3)=l$.
$$C=\left(\matrix{ 1&0&0\cr 0&0&1\cr 0&2l-1&0\cr}\right)\Longrightarrow
{\rm (D}_{l\ne\frac{1}{2}}-{\rm S}-2)= \left(\matrix{ 0&0&(2l-1)e_1\cr
0&e_1&0\cr  2le_1
&le_2&(2l-1)\cr}\right).$$

\noindent Case (BII):\quad (BII-1)\quad $f_{33}(e_3)=\lambda,\lambda\ne 0$.
$$\forall \lambda\ne 0,\;\;C=\left(\matrix{ 0&1&0\cr 1&0&0\cr
0&0&\lambda\cr}\right)\Longrightarrow
{\rm (D}_l-11)_\lambda = \left(\matrix{ 0&0&\lambda e_1\cr 0&0&\lambda
e_2\cr  (\lambda+1)e_1
&(\lambda+l)e_2&\lambda e_3\cr}\right),\lambda \ne 0;$$
when $\lambda=-1$, there is an additional equivalent class:
$$C=\left(\matrix{ 0&1&0\cr 1&0&0\cr 0&1&-1\cr}\right)\Longrightarrow
{\rm (D}_l-12) = \left(\matrix{ 0&0&-e_1\cr 0&0&-e_2\cr  0
&(l-1)e_2&-e_3+e_1\cr}\right);$$
when $\lambda=-l$, there is an additional equivalent class:
$$C=\left(\matrix{ 0&1&0\cr 1&0&0\cr
\frac{1}{l}&0&-l\cr}\right)\Longrightarrow
{\rm (D}_l-13) = \left(\matrix{ 0&0&-le_1\cr 0&0&-le_2\cr
(1-l)e_1&0&-le_3+e_2\cr}\right);$$
when $\lambda=l-1$, there is an additional equivalent class:
$$C=\left(\matrix{ 0&1&0\cr 1&\frac{1}{1-l}&0\cr
0&0&l-1\cr}\right)\Longrightarrow
{\rm (D}_l-14) = \left(\matrix{ 0&0&(l-1)e_1\cr 0&0&(l-1)e_2+e_1\cr  \
le_1& (2l-1)e_2+e_1&(l-1)e_3\cr}\right);$$
when $\lambda=1-l$, there is an additional equivalent class:
$$C=\left(\matrix{ \frac{1}{l-1}&1&0\cr 1&0&0\cr
0&0&1-l\cr}\right)\Longrightarrow
{\rm (D}_l-15) = \left(\matrix{ 0&0&(1-l)e_1+e_2\cr 0&0&(1-l)e_2\cr
(2-l)e_1+e_2&e_2&(1-l)e_3\cr}\right).$$

(BII-2)$_{\frac{1}{2}}$:\quad
(only for $l=\frac{1}{2}$)\quad $f_{33}(e_3)=-\frac{1}{2}$. Besides ${\rm
(D_{\frac{1}{2}}-11)}_{-\frac{1}{2}}$,
${\rm (D_{\frac{1}{2}}-13)}$ and ${\rm (D_{\frac{1}{2}}-14)}$, there is an
additional equivalent class:
$$C=\left(\matrix{ 0&1&0\cr 1&2&0\cr
2&0&-\frac{1}{2}\cr}\right)\Longrightarrow
{\rm (D_{\frac{1}{2}}-S-4)}= \left(\matrix{ 0&0&-\frac{1}{2}e_1\cr
0&0&-\frac{1}{2}e_2+e_1\cr  \frac{1}{2}e_1
&e_1&-\frac{1}{2}e_3+e_2\cr}\right).$$

\noindent Case ${\rm (BIII)}_{\frac{1}{2}}$: $f_{33}(e_3)=\lambda,\lambda\ne
0$.
$$\forall \lambda\ne 0,\;\;C=\left(\matrix{ 1&0&0\cr 0&1&0\cr
0&0&\lambda\cr}\right)\Longrightarrow
{\rm (D_{\frac{1}{2}}-S-5)}_\lambda = \left(\matrix{ 0&0&\lambda e_1\cr
0&e_1&\lambda e_2\cr  (\lambda+1)e_1
&(\lambda+\frac{1}{2})e_2&\lambda e_3\cr}\right),\lambda \ne 0;$$
when $\lambda=-1$, there is an additional equivalent class:
$$C=\left(\matrix{ 1&0&0\cr 0&1&0\cr 1&0&-1\cr}\right)\Longrightarrow
{\rm (D_{\frac{1}{2}}-S-6)} = \left(\matrix{ 0&0&-e_1\cr 0&e_1&-e_2\cr  0
&- \frac{1}{2}e_2&-e_3+e_1\cr}\right).$$
when $\lambda=\frac{1}{2}$, there is an additional equivalent class:
$$C=\left(\matrix{ 0&1&0\cr 0&0&1\cr 1&0&0\cr}\right)\Longrightarrow
{\rm (D_{\frac{1}{2}}-S-7)} = \left(\matrix{ 0&e_3&0\cr e_3&e_1&0\cr  e_1
&\frac{1}{2}e_2&\frac{3}{2}e_3\cr}\right).$$

\noindent Case (BIV):\quad (BIV-1)\quad $(l\ne -1)$\quad
$f_{22}(e_3)=\lambda,
f_{33}(e_3)=l, \lambda\ne 0$.
$$\forall \lambda\ne 0, \;\;C=\left(\matrix{ 1&0&0\cr 0&0&1\cr
0&\lambda&0\cr}\right)\Longrightarrow
{\rm (D}_{l\ne -1}-16)_\lambda = \left(\matrix{ 0&0&\lambda e_1\cr 0&0&0\cr
(\lambda+1) e_1&le_2&\lambda e_3\cr}\right),\lambda\ne 0;$$
when $\lambda=-1$, there is an additional equivalent class:
$$C=\left(\matrix{ 1&0&0\cr 0&0&1\cr 1&-1&0\cr}\right)
\Longrightarrow{\rm (D}_{l\ne -1}-17) = \left(\matrix{ 0&0&-e_1\cr 0&0&0\cr
0
&le_2&-e_3+e_1\cr}\right).$$
when $\lambda=-1$, there is an additional equivalent class:
$$C=\left(\matrix{ 0&1&0\cr 0&0&1\cr 1&0&0\cr}\right)
\Longrightarrow{\rm (D}_{l\ne -1}-18) = \left(\matrix{ e_3&0&0\cr 0&0&0\cr
e_1
&le_2&2e_3\cr}\right).$$

(BIV-2)\quad $(l\ne -1)$\quad $f_{22}(e_3)=l-1, f_{33}(e_3)=1$.
$$C=\left(\matrix{ \frac{1}{l-1}&0&1\cr 1&0&0\cr 0&1&0\cr}\right)
\Longrightarrow{\rm (D}_{l\ne -1}-19) = \left(\matrix{ 0&0&e_2\cr 0&0&0\cr
e_1+e_2
&le_2&(l-1)e_3\cr}\right).$$

\noindent Case (BV):\quad ($l\ne -1$)\quad there does not exist any
bijective 1-cocycle;

\noindent Case (BVI):\quad ($l\ne -1$)\quad $f_{33}(e_3)=l$.
$$C=\left(\matrix{ 0&1&0\cr l^2&0&0\cr 0&0&l\cr}\right)\Longrightarrow
{\rm (D}_{l\ne -1}-20) = \left(\matrix{ 0&0&le_1\cr 0&0&0\cr (1+l) e_1
&le_2&le_3+e_2\cr}\right).$$

\noindent Case (CI):\quad ($l\ne -1$)\quad $f_{33}(e_3)=l-2$.
$$C=\left(\matrix{ 0&1&0\cr 1&0&0\cr 0&0&l-2\cr}\right)\Longrightarrow
{\rm (D}_{l\ne -1}-21)= \left(\matrix{ e_2&0&(l-2)e_1\cr 0&0&0\cr (l-1) e_1
&le_2&(l-2)e_3\cr}\right).$$

{\bf Remark 1}\quad It is easy to see that we can extend the extent
of some parameters appearing in above left-symmetric algebras:
$$({\rm D}_l-1)_{\lambda=0}\cong ({\rm D}_l-4)_{\lambda=0}\cong
({\rm D}_l-11)_{\lambda=0};\;\;({\rm
D_{\frac{1}{2}}-S-1})_{\lambda=0}\cong({\rm
D_{\frac{1}{2}}-S-5})_{\lambda=0}.$$

{\bf Remark 2}\quad In some sense,  we can extend the value of $l$
to $l=1$ to get the classification of left-symmetric algebras on
${\cal D}_1$:
\begin{eqnarray*}
&&({\bar {\rm  D}}-1)_{\lambda}\cong ({\rm
D}_{l=1}-1)_{\lambda};\;\; ({\bar {\rm  D}}-2)\cong ({\rm
D}_{l=1}-2)_{\lambda}\cong ({\rm
D}_{l=1}-3);\\
&&({\bar {\rm  D}}-3)_{\lambda}\cong ({\rm
D}_{l=1}-4)_{\lambda}\cong ({\rm D}_{l=1}-16)_{\lambda};\;\;
({\bar {\rm  D}}-4)\cong ({\rm D}_{l=1}-5)\cong ({\rm D}_{l=1}-17);\\
&&({\bar {\rm  D}}-6)\cong ({\rm D}_{l=1}-7)\cong ({\rm D}_{l=1}-14)
                     \cong ({\rm D}_{l=1}-15)\cong({\rm D}_{l=1}-19);\\
&&({\bar {\rm  D}}-5)\cong ({\rm D}_{l=1}-6)\cong ({\rm
D}_{l=1}-18);\;\;
({\bar {\rm D}}-7)\cong ({\rm D}_{l=1}-8)\cong({\rm D}_{l=1}-20);\\
&& ({\bar {\rm D}}-8)\cong ({\rm D}_{l=1}-{\rm S}-1)\cong({\rm
D}_{l=1}-21);\;\;
({\bar {\rm D}}-9)\cong ({\rm D}_{l=1}-9)\cong  ({\rm D}_{l=1}-{\rm S}-2);\\
&&({\bar {\rm  D}}-10)\cong ({\rm D}_{l=1}-10);\;\;
  ({\bar {\rm  D}}-11)_\lambda \cong ({\rm D}_{l=1}-11)_\lambda;\;\;
({\bar {\rm  D}}-12)\cong ({\rm D}_{l=1}-12)\cong({\rm D}_{l=1}-13).
\end{eqnarray*}
However, for some cases, the corresponding bijective is quite
different. For example, $({\bar {\rm D}}-6)$ belongs to the case
(AV) of ${\cal D}_1$, but $({\rm D}_l-7)$ belongs to the case (AIV).

{\bf Remark 3}\quad Similarly, we can extend the value of $l$ to
$l=0$ to get certain left-symmetric algebras on ${\cal N}$:
\begin{eqnarray*}
&&({\rm  N}-1)_{\lambda}\cong ({\rm D}_{l=0}-1)_{\lambda};\;\; ({\rm
N}-5)\cong ({\rm D}_{l=0}-2)\cong ({\rm D}_{l=0}-13)\cong ({\rm
D}_{l=0}-20);\\
&& ({\rm  N}-6)\cong ({\rm D}_{l=0}-3);\;\; ({\rm
N}-9)_{\lambda}\cong ({\rm D}_{l=0}-4)_{\lambda};\;\; ({\rm
N}-5)\cong ({\rm D}_{l=0}-5);\;\;
({\rm  N}-10)\cong ({\rm D}_{l=0}-6);\\
&&({\rm  N}-11)\cong ({\rm D}_{l=0}-7);\;\; ({\rm  N}-15)\cong ({\rm
D}_{l=0}-8);\;\;
({\rm  N}-16)\cong ({\rm D}_{l=0}-{\rm S}-1);\\
&&({\rm  N}-36)\cong ({\rm D}_{l=0}-9);\;\; ({\rm  N}-35)\cong ({\rm
D}_{l=0}-10);\;\;
({\rm  N}-33)\cong ({\rm D}_{l=0}-{\rm S}-2);\\
&&({\rm  N}-38)_{\lambda}\cong ({\rm D}_{l=0}-11)_{\lambda};\;\;
({\rm  N}-39)\cong ({\rm D}_{l=0}-12);\;\;
({\rm  N}-40)\cong ({\rm D}_{l=0}-14);\\
&&({\rm  N}-42)\cong ({\rm D}_{l=0}-15);\;\; ({\rm  N}-18)_\lambda
\cong ({\rm D}_{l=0}-16)_\lambda;\;\;
({\rm  N}-19)\cong ({\rm D}_{l=0}-17);\\
&&({\rm N}-20)\cong ({\rm D}_{l=0}-18);\;\; ({\rm  N}-21)\cong ({\rm
D}_{l=0}-19);\;\; ({\rm  N}-44)\cong ({\rm D}_{l=0}-21).
\end{eqnarray*}
All above algebras satisfy the condition:
$f_{11}(e_1)=f_{22}(e_1)=f_{33}(e_1)=0$. However, there are certain
left-symmetric algebras on ${\cal N}$ satisfying this condition such
as (N-12), (N-17), (N-34) and (N-41) which cannot be obtained from
${\cal D}_l$ as $l=0$.

{\bf Proposition 3.14}\quad With the notations as above, among the
left-symmetric algebras on ${\cal D}_l$, we have

a) Associative algebras: $({\rm D}_{-1}-4)_1$;

b) Transitive left-symmetric algebras:$({\rm D}_{l}-1)_0$;$({\rm
D_{-1}-10})$; $({\rm D_{\frac{1}{2}}-S}-1)_0$;

c) Novikov algebras:$({\rm D}_{l}-1)_0$;
$({\rm D}_{l}-11)_\lambda$ ($\lambda\ne 0$);$({\rm D_l-12})$; $({\rm
D_l-13})$; $({\rm D_{\frac{1}{2}}-14})$;
$({\rm D_{\frac{1}{2}}-S}-1)_0$;
$({\rm D_{\frac{1}{2}}-S}-4)$;
$({\rm D_{\frac{1}{2}}-S}-5)_\lambda$ ($\lambda\ne0$);$({\rm
D_{\frac{1}{2}}-S}-6)$;

d) Bi-symmetric algebras:$({\rm D}_{-1}-4)_1$;$({\rm D}_{-1}-5)$;$({\rm
D}_{-1}-8)$;
$({\rm D_{-1}-T}-1)$;

e) Simple left-symmetric algebras on:$({\rm D}_{l}-10)$;$({\rm
D_{\frac{1}{2}}-S}-7)$.

\subsection{The left-symmetric algebras on ${\cal E}$}

The automorphism group of ${\cal E}$ is
$${\rm Aut}({\cal E})=\{
\left( \matrix{ a_{11} & 0 & 0\cr  a_{21} & a_{11} & 0\cr a_{31} &
a_{32}& 1\cr}\right)|a_{11}^2\ne 0\}.\eqno (3.8)$$

{\bf Proposition 3.15}\quad The equivalent classes of the
representations of ${\cal E}$ are divided into the following cases :

\noindent ${\rm (AI)}\; f(e_1)=0, f(e_2)=0, f(e_3)=\left(\matrix{
f_{11}(e_3)&
0 &0\cr 0 &f_{22}(e_3) &0\cr
0 & 0&f_{33}(e_3)\cr}\right).$

\noindent ${\rm (AII)}\; f(e_1)=0, f(e_2)=0, f(e_3)=\left(\matrix{
f_{11}(e_3)&
0 &0\cr 1&f_{11}(e_3) &0\cr
0 & 0&f_{33}(e_3)\cr}\right).$

\noindent ${\rm (AIII)}\; f(e_1)=0, f(e_2)=0, f(e_3)=\left(\matrix{
f_{11}(e_3)&
0 &0\cr 1&f_{11}(e_3) &0\cr
0 & 1&f_{11}(e_3)\cr}\right).$

\noindent ${\rm (AIV)}\; f(e_1)=0,
f(e_2)=\left(\matrix{ 0 &0&0\cr 1&0&0\cr 0&0&0\cr}\right),
f(e_3)=\left(\matrix{ f_{22}(e_3)+1&
0 &0\cr 0 &f_{22}(e_3) &0\cr
0 & 0&f_{33}(e_3)\cr}\right).$

\noindent ${\rm (AV)}\; f(e_1)=0,
f(e_2)=\left(\matrix{ 0 &0&0\cr 1&0&0\cr 0&0&0\cr}\right),
f(e_3)=\left(\matrix{ f_{33}(e_3)&
0 &0\cr 0 &f_{33}(e_3)-1 &0\cr
1 & 0&f_{33}(e_3)\cr}\right).$

\noindent ${\rm (AVI)}\; f(e_1)=0,
f(e_2)=\left(\matrix{ 0 &0&0\cr 0&0&0\cr 0&1&0\cr}\right),
f(e_3)=\left(\matrix{ f_{33}(e_3)&
0 &0\cr 0 &f_{33}(e_3)+1 &0\cr
1 & 0&f_{33}(e_3)\cr}\right).$

\noindent ${\rm (AVII)}\; f(e_1)=0,
f(e_2)=\left(\matrix{ 0 &0&0\cr 1&0&0\cr 0&1&0\cr}\right),
f(e_3)=\left(\matrix{ f_{33}(e_3)+2&
0 &0\cr 0 &f_{33}(e_3)+1 &0\cr
0 & 0&f_{33}(e_3)\cr}\right).$

\noindent ${\rm (BI)}\; f(e_1)=\left(\matrix{ 0 &0&0\cr 0&0&0\cr
1&0&0\cr}\right),
f(e_2)=\left(\matrix{ 0 &0&0\cr 0&0&0\cr 0&1&0\cr}\right),
f(e_3)=\left(\matrix{ f_{33}(e_3)+1&
0 &0\cr 1 &f_{33}(e_3)+1 &0\cr
0 & 0&f_{33}(e_3)\cr}\right).$

\noindent ${\rm (BII)}\; f(e_1)=\left(\matrix{ 0 &0&0\cr 0&0&0\cr
1&0&0\cr}\right),
f(e_2)=\left(\matrix{ 0 &0&0\cr 1&0&0\cr 0&0&0\cr}\right),
f(e_3)=\left(\matrix{ f_{33}(e_3)+1&
0 &0\cr 0 &f_{33}(e_3) &0\cr
0 & -1&f_{33}(e_3)\cr}\right).$

{\bf Proposition 3.16}\quad For a representation of ${\cal E}$ given
in the above cases respectively, there exist bijective 1-cocycles if
and only if it is equivalent to one of the following corresponding
cases:

\noindent Case (AI): there does not exist any bijective 1-cocycle;

\noindent Case (AII):\quad $f_{11}(e_3)=1,f_{33}(e_3)=\lambda,\lambda\in
{\bf C}$.
$$\forall \lambda\in {\bf C},\;\;C=\left(\matrix{ 1&0&0\cr 0&1&0\cr
0&0&1\cr}\right)\Longrightarrow
{\rm (E-1)}_\lambda = \left(\matrix{ 0&0&0\cr 0&0&0\cr  e_1
&e_1+e_2&\lambda e_3\cr}\right),\lambda\in {\bf C}.$$

\noindent Case (AIII):\quad $f_{11}(e_3)=1$.
$$C=\left(\matrix{ 1&0&0\cr 0&1&0\cr 0&0&1\cr}\right)\Longrightarrow
{\rm (E-2)}= \left(\matrix{ 0&0&0\cr 0&0&0\cr  e_1
&e_1+e_2&e_3+e_2\cr}\right).$$

\noindent Case (AIV):\quad $f_{22}(e_3)=0, f_{33}(e_3)=1$.
$$C=\left(\matrix{ 1&0&0\cr 0&0&1\cr 0&-1&0\cr}\right)\Longrightarrow
{\rm (E-3)}= \left(\matrix{ 0&0&0\cr 0&0&-e_1\cr  e_1
&e_2&0\cr}\right).$$

\noindent Case (AV):\quad $f_{33}(e_3)=1$.
$$  C_\lambda=\left(\matrix{ \frac{1}{\lambda}&0&0\cr
0&0&\frac{\lambda+1}{\lambda}\cr 0&1&0\cr}\right), \lambda\ne -1, 0\quad
\Longrightarrow{\rm (E-4)}_\lambda= \left(\matrix{ 0&0&0\cr 0&0&\lambda
e_1\cr  e_1
&e_2+(\lambda+1)e_1&0\cr}\right),\lambda\ne 0,-1.$$

\noindent Case (AVI):\quad $f_{11}(e_3)=1$.
$$C=\left(\matrix{ 1&0&0\cr 0&1&1\cr 0&0&1\cr}\right)\Longrightarrow
{\rm (E-5)}= \left(\matrix{ 0&0&0\cr 0&e_2-e_3&e_2-e_3\cr  e_1
&e_1+2e_2-e_3&e_3+e_1\cr}\right)\cong\left(\matrix{ 0&0&0\cr 0 &e_3&0\cr
e_1& e_1+e_2& 2e_3\cr}\right).$$

\noindent Case (AVII):\quad $f_{33}(e_3)=-1$.
$$C=\left(\matrix{ 1&0&0\cr 0&1&0\cr 0&-1&-1\cr}\right)\Longrightarrow
{\rm (E-6)}= \left(\matrix{ 0&0&0\cr 0&e_1&-e_1-e_2\cr  e_1
&0&-e_3-e_2\cr}\right).$$

\noindent Case (BI):\quad $f_{33}(e_3)=\lambda,\lambda\ne 0$.
$$\forall \lambda\ne 0,\;\;C=\left(\matrix{ 1&0&0\cr 0&1&0\cr 0&0&\lambda
\cr}\right)\Longrightarrow
{\rm (E-7)}_\lambda = \left(\matrix{ 0&0&\lambda e_1\cr 0&0&\lambda e_2\cr
(\lambda+1) e_1&e_1+(\lambda+1)e_2 &\lambda e_3\cr}\right),\lambda\ne 0;$$
when $\lambda =-1$, there is an additional equivalent class:
$$C=\left(\matrix{ 1&0&0\cr 0&1&0\cr -1&1&-1 \cr}\right)\Longrightarrow
{\rm (E-8)}= \left(\matrix{ 0&0&-e_1\cr 0&0&-e_2\cr
0&e_1&-e_3+e_2\cr}\right).$$

\noindent Case (BII):\quad $f_{11}(e_3)=1$.
$$C=\left(\matrix{ 1&0&0\cr 1&1&0\cr -1&0&1\cr}\right)\Longrightarrow
{\rm (E-9)}= \left(\matrix{ 0&0&e_1\cr 0&e_1&0\cr  2e_1
&e_1+e_2&e_3-e_2\cr}\right).$$

{\bf Remark }\quad It is easy to see that we can extend the extent
of some parameters appearing in above left-symmetric algebras:
$$({\rm E-1})_{\lambda=0}\cong ({\rm E-4})_{\lambda=0}\cong
({\rm E-7})_{\lambda=0};\;\;({\rm E-4})_{\lambda=-1}\cong({\rm
E-3}).$$

{\bf Proposition 3.17}\quad With the notations as above. Among the
left-symmetric algebras on ${\cal E}$, we have

a) There is not any associative algebra on ${\cal E}$;

b) Transitive left-symmetric algebras:(E-1)$_0$;(E-3); (E-4)$_{\lambda}$
($\lambda\ne 0,-1$);

c) Novikov algebras: (E-1)$_0$; (E-7)$_\lambda$($\lambda\ne 0$); (E-8);

d) There is not any bi-symmetric algebra on ${\cal E}$;

e) There is not any simple left-symmetric algebra on ${\cal E}$.

\section*{Acknowledgements}

The author thanks Professors D. Burde, I.M. Gel'fand, B.A.
Kupershmidt and P. Etingof for useful suggestion and great
encouragement. In particular, the author is grateful for Professor
P. Etingof who invited the author to visit MIT and the work was
begun there. The author also thanks Professors J. Lepowsky, Y.Z.
Huang and H.S. Li for the hospitality extended to him during his
stay at Rutgers, The State University of New Jersey and for valuable
discussions. This work was supported in part by S.S. Chern
Foundation for Mathematical Research, the National Natural Science
Foundation of China (10571091, 10621101), NKBRPC (2006CB805905),
Program for New Century Excellent Talents in University and K.C.
Wong Education Foundation.

\end{document}